\documentclass[a4paper,11pt]{article}

\usepackage{graphicx}
\usepackage{bm}

\usepackage{authblk}
\usepackage{amscd}
\usepackage{amsfonts}
\usepackage[numbers]{natbib}
\usepackage{amsthm}
\usepackage{amssymb}
\usepackage{amsmath}
\usepackage{algpseudocode}
\usepackage{algorithm}
\usepackage{cancel}
\usepackage{chngcntr}
\usepackage{caption}
\usepackage{enumerate}
\usepackage{epsfig}
\usepackage{epic}
\usepackage{eepic}
\usepackage{enumitem}
\usepackage{etoolbox}
\usepackage{fancyhdr}
\usepackage{fullpage}
\usepackage[a4paper, total={6in, 8in}]{geometry}
\usepackage{graphicx}
\usepackage{here}
\usepackage[utf8]{inputenc}
\usepackage{lscape}
\usepackage{listings}
\usepackage{longtable}
\usepackage{natbib}
\usepackage{nccmath}
\usepackage{parskip}
\usepackage{rotating}
\usepackage{subfigure}
\usepackage{tabularx}
\usepackage{textcomp}
\usepackage{tikz}
\usepackage{titlesec}
\usepackage{threeparttable}
\usepackage{varwidth}
\usepackage{framed,multirow}

\usepackage{xcolor}

\newcommand{\vx}{\vec{x}}
\newcommand{\tpar}{\tilde{\partial}}

\newcommand{\pad}[2]{\frac{\partial #1}{\partial #2}}
\newcommand{\vepst}{\tilde{\varepsilon}}
\newcommand{\veps}{\varepsilon}

\title{Multiscale Modeling of Sorption Kinetics}

\author[1]{Clarissa Astuto}
\author[2]{Antonio Raudino}
\author[3]{Giovanni Russo}

\affil[1]{King Abdullah University of Science and Technology (KAUST), 4700, Thuwal, Saudi Arabia}
\affil[2]{Department of Chemical Sciences,  University of Catania, Italy}
\affil[3]{Department of Mathematics and Computer Science, University of Catania, Italy}

\begin{document}
\maketitle

\begin{abstract}
In this paper we propose and validate a multiscale model for the description of particle diffusion in presence of trapping boundaries. We start from a drift-diffusion equation in which the drift term describes the effect of bubble traps, and is modeled by a short range potential with an attractive term and a repulsive core. 
The interaction of the particles attracted by the bubble surface is simulated by the Lennard-Jones potential that simplifies the capture due to the hydrophobic properties of the ions.
In our model the effect of the potential is replaced by a suitable boundary condition derived by mass conservation and asymptotic analysis. The potential is assumed to have a range of small size $\varepsilon$. An asymptotic expansion in the $\varepsilon$ is considered, and the boundary conditions are obtained by retaining the lowest order terms in the expansion. 

Another aspect we investigate is saturation effect coming from high concentrations in the proximity of the bubble surface. The validity of the model is carefully checked with several tests in 1D, 2D and different geometries. 
\end{abstract}

	\section{Introduction}
\label{sec:introduction}
The diffusion of particles in presence of (moving or static) traps is an interesting topic because of its applications to different fields such as chemistry, physics and biology. An interesting case is to consider the trapping of diffusing particles near the surface of the traps. The problem of surfactant diffusion, that we consider in this paper, has been investigated by several authors \cite{BERG1977193,DIAMANT2001259,MILLER2017115,MILLER1994249,WIEGEL1983283} and in each work the transport phenomena in ionic solutions requires the calculation of the space and time evolution of the concentration of negative (anions) and positive (cations) species diffusing in a confined space at a time $t$, as we did in \cite{CiCP-31-707}. 

In this work we study the behaviour of one species at a time because here we focus on the space multiscale aspects of the problem, while many authors have focused on the time multiscale challenge of the phenomenon \cite{MILLER1994249,BEVERUNG199959}. The most frequently used strategy to describe the dynamic properties of liquid adsorption layers is the dynamic surface or the interfacial tension. Many different techniques exist to cover the time interval from milliseconds to hours 
and days.  

Introducing the local concentration of ions $c= c(\vec{x},t)$, its time evolution in a fluid  is governed by the conservation law
\begin{equation}
	\displaystyle \frac{\partial c\left(\vec{x},t\right)}{\partial t}=-\nabla\cdot J(\vec{x},t).
	\label{equation_flux}
\end{equation}
In absence of bubbles, the flux $J$ represents the standard diffusion  {given by Fick's law}  $J = - D \nabla c$, where $D$ is the diffusion coefficient. In the presence of bubbles (traps), the flux expression is augmented to include a {\it drift} effect
\begin{eqnarray}
	\displaystyle J=\ -D\left(\nabla c +\ \frac{1}{k_BT}c\nabla V\right)   	
	\label{equation_definition_flux}	
\end{eqnarray}
where $k_B$ is the Boltzmann's constant, $T$ is the absolute temperature (assumed to be constant) and $V(x)$ is a suitable potential function that models the \textit{attractive-repulsive} behavior of the bubble with the particles. In particular,  {when a particle is near the bubble, it is attracted towards its surface, and when it is at a very short distance from the surface of the bubble, the potential is designed to repulse the particle in order to simulate impermeability, therefore} trapping a portion of particles in the proximity of the surface of the bubble.

The potential $V$, in presence of positive and negative diffusing particles, may depend also on $c$ to describe the electrostatic interaction between the two species for high concentrations. In \cite{CiCP-31-707} we partitioned the potential as: 
\begin{eqnarray}
	V =\ V_{\rm{ion-bubble}} +\ V_{\rm{ion-ion}}  
	\label{equation_definition_ext}
\end{eqnarray}
where the term $V_{\rm{ion-bubble}}$ describes the interaction between the bubble with the ions, located at a generic position $r$ (Fig.~\ref{fig:bubble}), and the (possibly oscillating) interface, while the term $V_{\rm{ion-ion}}$ accounts for the interaction among the diffusing ions. {In this work we focus on modeling only the interaction of anions with the bubble, so we do not take into account the term $V_{\rm{ion-ion}}$, describing the interaction with the cations. Such term can indeed be neglected if only anions are present in the fluid, or, more realistically, when the Debye length associated to the two types of ions is much smaller than the size of the bubble (quasi-neutral limit, see \cite{CiCP-31-707}). The more realistic case of bipolar transport, which includes the Coulomb interaction between ions, will be subject of future investigation.}


{In this work we also include the study of the \emph{saturation} effect  of the system} when high concentrations get stuck to the surface of the trap (\cite{1946JChPh,Liu2000}). Standard expressions for the flux satisfies a Fermi-Dirac form for the mixing entropy: \cite{chavanis,Martzel_2001,Ribeiro}
$S_{\rm mix} = -k(c\log(c) + (1-c)\log(1-c))$  {providing the following expression for the flux}
\begin{eqnarray}
	\displaystyle J=\ -D\left(\nabla c +\ \frac{1}{k_BT}c(1-c)\nabla V\right)  
	\label{equation_nonlinear_flux}
\end{eqnarray}
and it leads to a non-linear equation in $c$.  {Note that for low concentration this expression reduces to Eq.~\eqref{equation_definition_flux}.} The main reason why studies take into account saturation is that it seems to be 
{relevant near the trap, where local concentrations may become non negligible.} Several authors have demonstrated that it is more appropriate to
use models with nonlinear saturating mobility instead of the classical linear mobility models, both for a single  {species} (\cite{SIMPSON2009399,doi:10.1080/22054952.2009.11464027,Giacomin_1997,474717002cc1428ca1a16b2dbb8103ee}) and also for multiple species (\cite{SIMPSON2009399,Painter2009}). 
In \ref{section_saturation} we see the effect of the flux defined in Eq.~\eqref{equation_nonlinear_flux} that is different from the one acting in dilute solutions where one can set $1-c\approx 1$. 



{In the linear model there is no upper bound on the concentration $c$, while in the model which makes use of Eq.~\eqref{equation_nonlinear_flux} the natural bound $c<1$ is always satisfied. }

Numerical results and the resulting saturation curve show how the concentration can not go over certain values ($c\leq 1$ in our case).  This effect is known to exist when the ions are exposed to a sufficiently strong potential (trap) but it can not be obtained in a limit of \emph{quasi-neutrality} \cite{CiCP-31-707}.	

Let us briefly summarize similar models dealing with the same issue that have been widely discussed in the literature, because of numerous applications in different problems in chemistry, physics and biology. For the sake of simplicity, people consider a semi-infinite mono-dimensional system, i.e.: $0 < x < \infty$.

Since the potential acting on the diffusing particles is mainly located in a narrow region near the left boundary and null otherwise (see Fig.~\ref{figure_potential_V_1D}), the transport equation can be modeled by a standard diffusion equation (Eq.~\eqref{eq_1d_inf}) together with the condition that the flux of particles $\displaystyle - J = D \left. \frac{ \partial c}{\partial x}\right|_{x=0} $, entering the left boundary, equals the time variation of surface concentration ($-J = \partial \Gamma / \partial t$) at the interface (continuity equation). So doing, the effect of the potential is adsorbed into a boundary condition at $x = 0$:
\begin{eqnarray}
	\label{eq_1d_inf}
	\displaystyle \frac{\partial c }{\partial t} &=& D \frac{\partial^2 c}{\partial x^2} \quad {\rm for} \quad x\in[0,+\infty[, t>0\\
	\label{eq_1d_inf2}
	\displaystyle \frac{\partial \Gamma}{ \partial t} &=& D\left.\frac{\partial c}{\partial x}\right|_{x=0} \quad {\rm for} \quad x=0, t>0
\end{eqnarray}
and initial condition
\begin{equation}
	c(x,0) = c_0 \quad x>0, t=0 
\end{equation}
Eq.(\ref{eq_1d_inf}-\ref{eq_1d_inf2})  must be supplemented by the further conditions: 
\begin{eqnarray}
	c|_{x=0} &=& f(\Gamma) \label{eq:f(Gamma)}\\
	c|_{x=\infty} &=& c_0
\end{eqnarray}
In Eq.~\eqref{eq_1d_inf2}, $f(\Gamma)$ is a prescribed function (see \cite{DIAMANT2001259,MILLER2017115,MILLER1994249,WIEGEL1983283,BEVERUNG199959}) (derived by empirical or theoretical arguments) relating $ c|_{x=0}$ (the diffusant concentration near the surface satisfying the diffusion equation \eqref{eq_1d_inf}) to  the still unknown diffusant concentration $\Gamma$ at the interface between the two adjacent phases, Eq.~\eqref{eq:f(Gamma)}. $\Gamma$ is an important quantity because it is related to experimentally accessible parameters like the surface tension. 
Other models in the literature (see, e.g.,\cite{item_1932442}), are based on a similar assumption because in addition to the bulk properties, the interface characteristics can affect drastically the motion of substances in the liquid. The presence of any little surfactant contamination can significantly affect the bubble motion and shape. Sorption effect is reversible: in some part of the surface particles may be absorbed, while in other part they may be desorbed, according to the sign of $\partial c/\partial x$. These adsorption and desorption fluxes at the bubble surface are governed by diffusion and lead to the formation of a so-called diffusion boundary layer adjacent to the surface \cite{item_1932442}.

In the simplest approximation (investigated by Sutherland a long time ago \cite{suth}) it is assumed: $f(\Gamma) \approx K\Gamma$. Such an approximation leads to rather cumbersome but exact analytical results (the same reference of  S.S. Dukhin as before) for the time evolution of $\Gamma$.

Our model represents an advancement over the classical approaches because:

A) Our approach allows the derivation of Sutherland's constant $K$ from first principles in the case of dilute solutions.

B) It can be generalized to concentrated solutions, providing a nonlinear relation between the surface concentration and the concentration near the trap.

C) A detailed validation of the derived model with the original drift-diffusion model described by a potential is performed.

{D) The multiscale model provides tremendous savings in computational complexity over the model based on fully resolved potential, even if the latter is solved adopting Adaptive Mesh Refinement, as illustrated in detail in Appendix~\ref{A:AMR}.}

The plan of the paper is the following: in Sec.~\ref{section_multiscale} we derive the multiscale model in 1D and 2D and then we extend the condition to 3D. In Sec.~\ref{section_validation_multiscale} we validate the model comparing its numerical results with those obtained with the original drift-diffusion model defined by Eqs.~(\ref{equation_flux}-\ref{equation_definition_flux}). Sec.~\ref{section_saturation} is devoted to the treatment of the saturation effect, which becomes relevant for non-negligible concentrations near the bubble surface. In the final Sec.~\ref{sec:conclusions} we draw some conclusions.

\section{Modeling multiscale problem}
\label{section_multiscale}
Modeling surfactant diffusion and trapping is very challenging because of the multiple scales involved in space. The range of the attractive-repulsive core of a trap is of the order of nanometers, a length that is prohibitively small for detailed simulations. 

To overcome such difficulties we propose a \emph{multi-scale model}. The ions are modeled by a volumetric concentration $c(x,t)$ which satisfies a simple diffusion equation in the domain of the liquid, coupled with a superficial one, which describes the carriers trapped on the surface. The two concentrations are related by a mass balance equation which acts as a suitable boundary condition for the concentration in the bulk. 

The multiscale model has been tested in various geometries and dimensions $d$. For $d=1$ the domain is $[-\varepsilon,1]$, the action of the trap is defined with an external potential that is non negligible in a few $\varepsilon$ and its minimum is at the origin. For $d = 2$, we first consider a natural extension of the 1D case with a trap described by a potential that is constant in $y$-direction. Secondly, we consider a  2D region inside a square, external to a circular disk. The arbitrary boundary of the bubble is described by a level set function. The set of equations is then solved by a finite difference scheme on a regular Cartesian mesh, which makes use of regular grid points (inside the domain $\Omega$) and ghost points (neighbors of regular grid points which are out of $\Omega$). The equations on the regular points are obtained by finite difference discretization of the Laplacian operator, while the equations on each ghost point are given by interpolating the numerical solution near the ghost point and imposing the boundary condition on the interpolant at the orthogonal projection of the point on $\partial \Omega$ (see \cite{COCO2013464,COCO2018299} and Sec.~\ref{section_validation_multiscale} for more details).

Upon the validation in the geometries described below, it is observed that when the range of the interaction decreases, the agreement of the multiscale model and the detailed numerical solution of the full drift diffusion model becomes better and better. 


\subsection{1D Multiscale model}
\label{section_1Dmodel}
\begin{figure}[H]
	\centering
	\begin{minipage}
		{.49\textwidth}
		\centering
		\includegraphics[width=\textwidth]{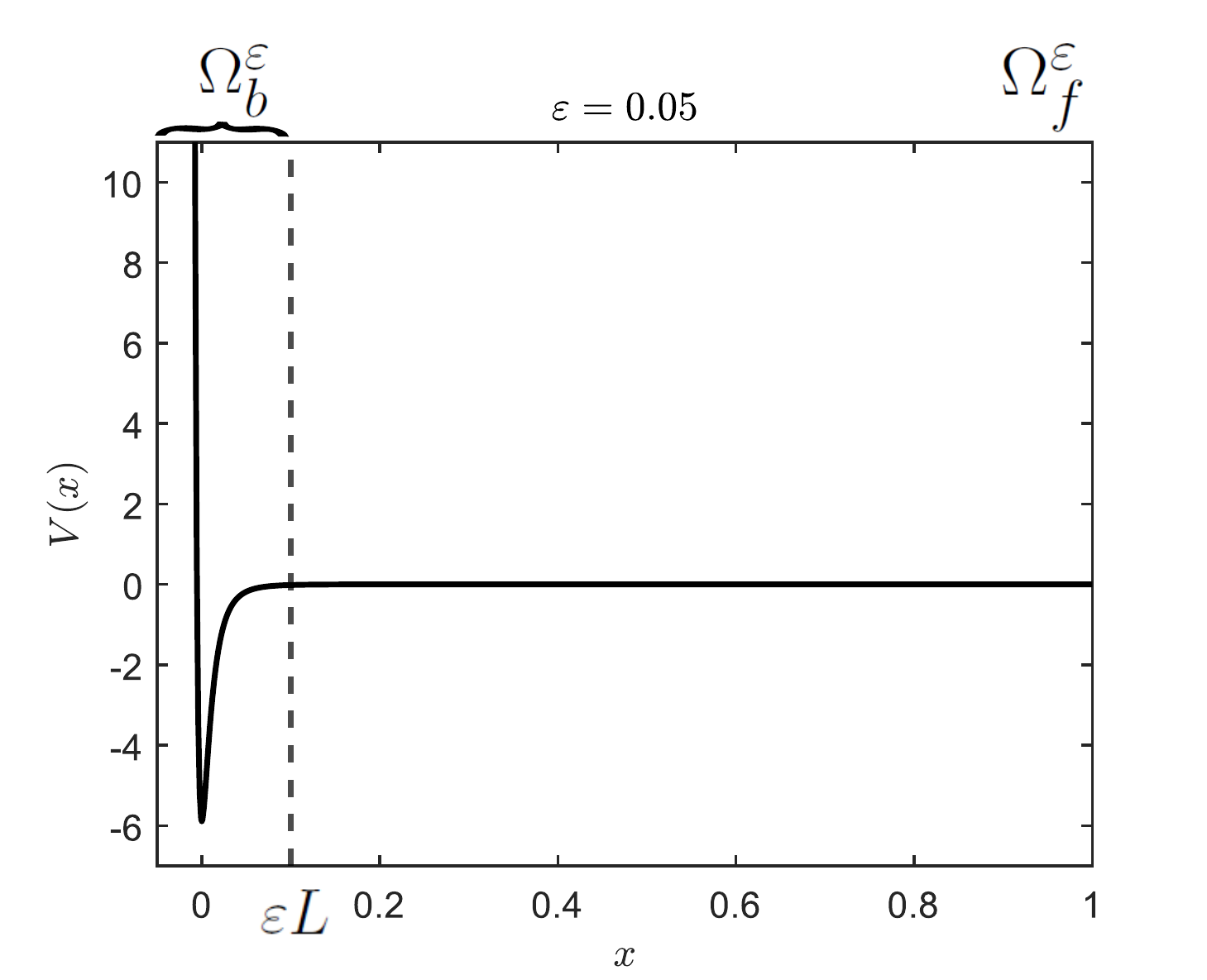}
	\end{minipage}
	\begin{minipage}
		{.49\textwidth}
		\includegraphics[width=\textwidth]{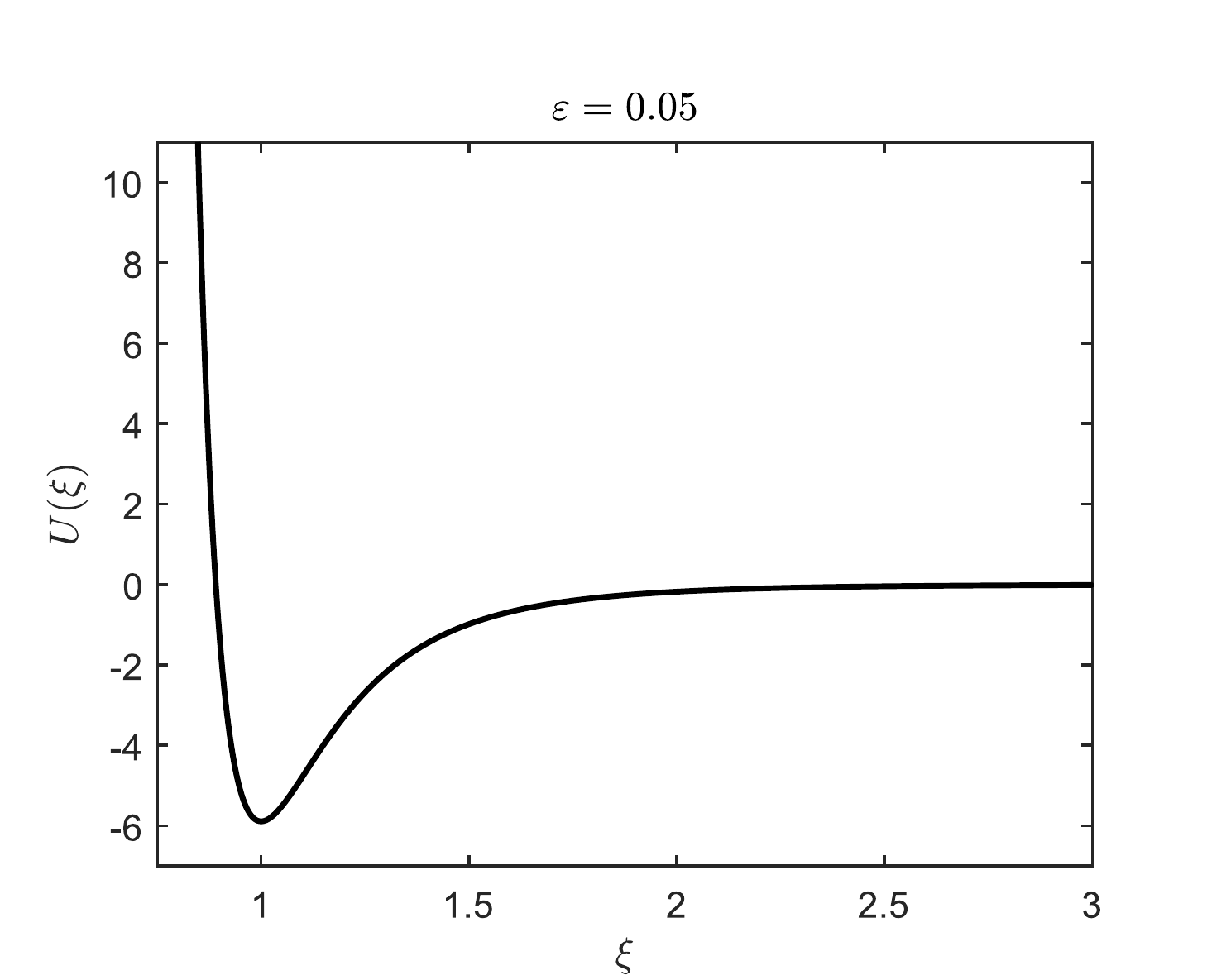}
	\end{minipage}	\caption{\textit{Representation of the external potentials $V(x)$, on the left, and $U(\xi)$, on the right, for $\varepsilon = 0.05$. $U(\xi)$ is obtained from a change of variable of $V(x)$, with $\xi \in [0,L+1]$. After $\xi = L+1$, the contribution of the potential $U(\xi)$ is negligible. For this reason, on the left the dashed line $x=\varepsilon L$ denotes the right boundary of $\Omega_b^\varepsilon$, the domain affected by the potential. Numerically, we notice an appropriate value for this quantity is $L=2$. }}
	\label{figure_potential_V_1D}
\end{figure}

{  In this subsection we derive the multiscale model in 1D space dimension. The drift-diffusion} Eqs.~\eqref{equation_flux}-\eqref{equation_definition_flux} in 1D read
\begin{eqnarray}
	\label{pde1d}
	\displaystyle \frac{\partial   c }{\partial t} + \frac{\partial  J}{ \partial x} = 0,
\end{eqnarray}
where
\begin{equation}
	\label{flux1d}
	\displaystyle J = - D\left( \frac{\partial  c  }{ \partial x} + \frac{1}{k_BT}  c  V'\right).
\end{equation}
We assume that {  the full domain is $\Omega = [-\varepsilon,1]$ and then we  partition it into a bubble domain $\Omega^\veps_b = [-\veps,L\veps]$, inside which the ions feel the effect of the potential, and} the fluid domain, which is not affected by the bubble, $\Omega^\varepsilon_f = [L\varepsilon,1],$  see Fig.~\ref{figure_potential_V_1D}. 
Therefore, the potential has effect only in $\Omega^\varepsilon_b$, while it is constant for $x \in \Omega^\varepsilon_f$, say $V(x)=0$. 

As a prototype potential we consider the Lennard-Jones potential (LJ), which describes attraction at long distances and repulsion at short distances due to Van der Waals and Pauli terms, respectively. A typical shape of the potential $V(x)$ is shown in Fig.~\ref{figure_potential_V_1D} 
and it takes the form
\begin{eqnarray*}
	V(x) = E\left( \left(\frac{x+\varepsilon}{\varepsilon}\right)^{-12} - 2\left(\frac{x+\varepsilon}{\varepsilon}\right)^{-6} \right).
\end{eqnarray*}
where $\varepsilon$ denotes the range of the potential and $E$ represents the depth of the well (see Fig.~\ref{figure_potential_V_1D}, left panel). It is convenient to adopt a non dimensional form of the potential, expressed as a function of the rescaled variable $\xi = 1 + x/\varepsilon \in [0,L+1]$:
\begin{equation}
	\label{expr_U_LJ}
	U(\xi) = \phi \left( \xi^{-12} - 2\xi^{-6} \right)
\end{equation}
where { $\displaystyle \phi = {E}/{k_BT}$ represents the ratio between the depth, $E$, of the potential well of the bubble and $k_B T$. Typical values of $\phi$ should range between 10 and 20 \cite{everett1976adsorption}. $L$ is a scaled distance beyond which the potential is negligible. 
}


Let us assume that the trap is located at $x=0$ and its ''thickness'' is of order of $\varepsilon$.
The original model requires to solve Eq.~\eqref{pde1d} in $\Omega^\varepsilon=\Omega^\varepsilon_b \cup \Omega^\varepsilon_f$ with initial condition $ c  =  c_0(x)$ when $t=0$ and with boundary conditions
\begin{equation}\label{bc1d}
	J=0 \quad {\rm at} \> x=-\varepsilon, {
		\rm and} \> x = 1
\end{equation}
to simulate wall effects. If the potential diverges at $x=-\varepsilon$ such a condition is not necessary, since the potential wall prevents the particles from approaching the boundary. This condition can be used in case we simplify the description of the potential by an attractive well and a vertical boundary.

The mathematical model proposed so far presents a multiscale nature, since the influence of the bubble on the particles is simulated in $\Omega^\varepsilon_b$ that is very small compared to the rest of the domain, leading to some extra computational effort to ensure that the correct behaviour is captured. 
In order to overcome this difficulty, we propose a multiscale model for $\varepsilon \ll 1$ in which the domain of the problem is reduced to $\Omega^\varepsilon_f$ and the effect of $V(x)$ on $\Omega^\varepsilon_b$ is approximated by a suitable boundary condition at the origin obtained as follows.
Let us assume that the potential $V_\varepsilon(x)$ depends on $\varepsilon$ but maintains the same functional form, i.e.~there exist two functions: ${U}\colon[0,L+1] \rightarrow \mathbb{R}$ (e.g. the LJ function defined in Eq.~\eqref{expr_U_LJ}) and $\mathcal{U}\colon\mathbb{R} \rightarrow \mathbb{R}$ that do not depend on $\varepsilon$ and such that
\begin{eqnarray}
	\displaystyle \mathcal{U}(\xi)=\left\{
	\begin{array}{l}
		+\infty \quad \xi\leq 0\\
		U(\xi) \quad \xi \in [0,L+1] \\
		0 \quad \quad \xi > L+1
	\end{array}
	\right.
\end{eqnarray}
It is clear that the solution $ c $ and the flux $J$ depend on $\varepsilon$ as well, and then the original problem can be written as
\begin{eqnarray}
	\label{eq_flux_eps}
	\displaystyle \frac{\partial  c _\varepsilon}{\partial t} + \frac{\partial J_\varepsilon}{\partial x} &=&0\quad \rm{in } \, \Omega^\varepsilon\\
	\displaystyle J_{\varepsilon} &=& -D\left(\frac{\partial  c _\varepsilon}{\partial x} + \frac{1}{k_B T} \,  c _\varepsilon  \, V'_\varepsilon \right)
\end{eqnarray}
For $x \in \Omega^\varepsilon_b$ we can use the scaled variable $\displaystyle\xi= 1 + \frac{x}{\varepsilon} \in [0,L+1]$ and then the flux in $\Omega^\varepsilon_b$ can be written as
\begin{equation} \label{newflux}
	J_\varepsilon = -D \frac{1}{\varepsilon}\left(\frac{\partial  c _\varepsilon}{\partial \xi}+ c _\varepsilon \, U' \right).  
\end{equation}
{Rewriting \eqref{eq_flux_eps} in the rescaled variable and using Eq.~\eqref{newflux} we obtain
	\begin{equation}
		\label{c_eps2}
		\frac{\partial  c_\varepsilon}{\partial t} = -\frac{1}{\veps}\pad{J_\veps}{\xi} = \frac{D}{\varepsilon^2}\frac{\partial}
		{\partial \xi} \left(\pad{c_\veps}{\xi} + c_\veps U'(\xi)\right).
	\end{equation}
}

The range of the scaled variable $\xi$ does not depend on $\varepsilon$ and then we can assume that $ c _\varepsilon(\xi,t)$ has the following expansion in $\Omega^\varepsilon_b$:
\begin{equation}\label{exprho}
	{c_\varepsilon(\xi,t) =  c ^{(0)}(\xi,t)+\veps^2  c^{(1)}(\xi,t)+
		\veps^4  c^{(2)}(\xi,t)+\cdots.}
\end{equation} 
{ 
	Inserting the expansion~\eqref{exprho} into Eq.~\eqref{newflux} we obtain the following expansion for the flux:
	\begin{equation}
		\veps J_\veps = J^{(0)} + \veps^2 J^{(1)}  + \veps^4 J^{(2)} + \cdots,
		\label{Jexpand}
	\end{equation}
	with 
	\[
	J^{(k)} = -D\left(\pad{ c^{(k)}}{\xi} + c^{(k)}U'(\xi) \right), \quad k\geq 0.
	\]
	Using expansion \eqref{Jexpand}  in \eqref{c_eps2} we obtain, to the various orders in $\veps$:
	\begin{eqnarray}
		O(\varepsilon^{-1})  :&  \phantom{\quad\frac{\partial  c^{(k)}}{\partial t} + } \pad{J^{(0)}}{\xi} & = 0 ,
		\label{eq:eps0}\\
		O(\varepsilon^{k})  : &  \quad\frac{\partial  c^{(k)}}{\partial t} + \pad{J^{(k)}}{\xi}
		&  =  0, \quad k\ge 0.\label{eq:epsk}
	\end{eqnarray}
	Eq.~\eqref{eq:eps0} states that  the lowest order flux $J^{(0)}$ is constant.\\
	Condition $\lim_{x\to-\varepsilon}J(x,t) = 0$, which is a consequence of the singularity of the potential, gives:
	\begin{equation}
		\frac{\partial  c ^{(0)}}{\partial \xi}+ U'(\xi) c ^{(0)} = 0. 
		\label{Jorder0} 
	\end{equation} 
}
Therefore, it is possible to integrate the equation
\begin{equation}
	\frac{1}{ c ^{(0)}}\frac{\partial  c^{(0)} }{\partial \xi} = - U'(\xi) 
\end{equation}
whose solution is
\begin{equation}
	c ^{(0)}(\xi,t) =  c ^{(0)}(L+1,t)\exp \left(-U(\xi) + U(L+1) \right) =  c ^{(0)}(L+1,t)\exp \left(-U(\xi)\right)
\end{equation}
since $U(L+1) = 0$.

\begin{figure}[H]
		\centering
			
		\includegraphics[width=0.7\textwidth]{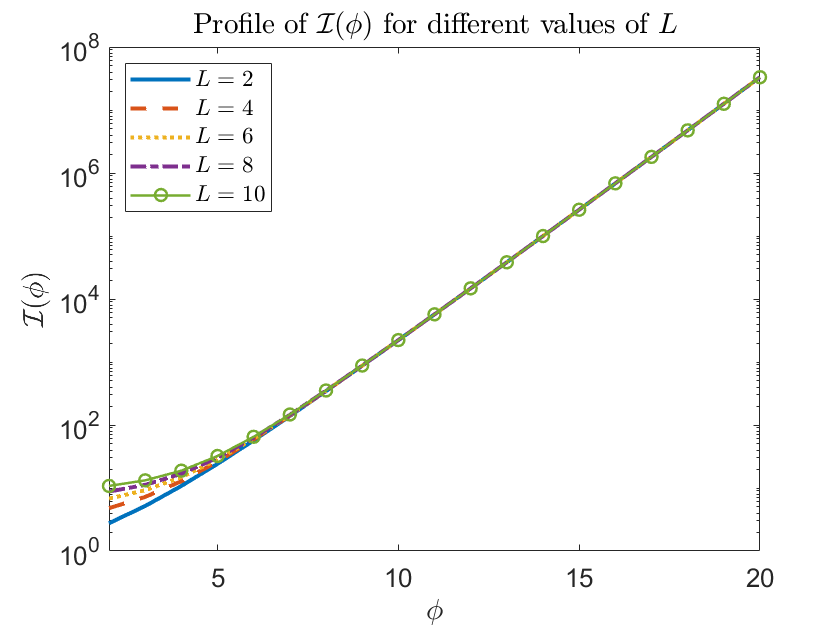}
		\caption{Plot of the profile of $\mathcal{I}(\phi)$, defined in Eq.~\eqref{eq_Iphi_mult}, {  with M = 3 and for different values of $L\in \{2,4,6,8,10\}$. This quantity increases with $\phi$ and we see that, for sufficiently large values of $\phi$, say $\phi>6$,  the values of $\mathcal{I}(\phi)$ are almost independent on $L$. }}
		\label{fig:Iphi_mult}
	\end{figure}

Now we assume that $V_x(\varepsilon L) = 0$, which is compatible {  with having} a smooth potential with compact support. Notice that, strictly speaking, LJ potential is not with compact support, however it is negligible small, together with its derivative, at distances which are {  sufficiently larger than its range $\varepsilon$ (see Table~\ref{table_L} and Fig.~\ref{fig:Iphi_mult}.)}

Integrating \eqref{pde1d} in the range of the potential, using the approximation $ c (x,t) =  c ^{(0)}(\xi,t), \, x \in \varepsilon[-1,L]$, the zero boundary condition \eqref{bc1d} at $x=-\varepsilon $ and that $V_\varepsilon(\varepsilon L)=V'_\varepsilon(\varepsilon L)=0$, we obtain
\[
\frac{d}{dt}\int_{- \varepsilon}^{\varepsilon L} c (x,t) \, dx +J_\varepsilon(\varepsilon L)  =  0,
\]
for which it follows
\begin{eqnarray}
	\nonumber
	\frac{\partial c (\varepsilon L,t)}{\partial t} \; \varepsilon \int_{0}^{L+1}\exp\left(-U(\xi) \right) d \xi - D \frac{\partial  c (\varepsilon L,t)}{\partial x} &=&0
\end{eqnarray}
that represents a boundary condition at $x=\varepsilon L$.
{Notice that, out of the range of the potential, the flux $J$ is given by just the diffusion term, thus Eq.~\eqref{eq:epsk} will not be used since for sufficiently small values of $\veps$ it is $c(x,t)\approx c^{(0)}(x,t)$, and Eq.~\eqref{eq:eps0} is used to provide an effective boundary condition for the diffusion equation.}

Summarizing, the \textit{multiscale model} can be stated as a limit model for $\varepsilon \to 0$ as follows:
\begin{eqnarray}\label{reduced1d}
	\frac{\partial  c }{\partial t} &=& D\frac{\partial^2  c }{\partial x^2} \quad {\rm {in} }\, x \in [0,1]\\
	\frac{\partial  c }{\partial x} &=& 0  \quad {\rm {at} }\, x = 1 \\ \label{BCt}
	M\frac{\partial  c }{\partial t} &=& D\frac{\partial  c }{\partial x}  \quad {\rm {at} }\, x = 0
\end{eqnarray}
where
\begin{equation}
	\label{expr_M}
	M=\varepsilon\int_{0}^{L+1}\exp\left(-U(\xi)\right)d\xi.
\end{equation}

{As $\veps\to 0$ one could let $L\to\infty$ in such a way that $L\veps\to 0$, so in the expression for $M$ the integral can be written from $0$ to $\infty$. 
	In practice we use a finite value for $L$. 
	Except where otherwise stated, we used $L=2$.  In Fig.~\ref{fig:Iphi_mult} we see that this is a convenient value for the parameter $L$, since the quantity $\mathcal{I}_L(\phi)$ approximately does not change increasing $L$. Also in Table~\ref{table_L} we see an agreement with the choice of the parameter $L$. We show that the action of the potential $U(\xi)$ is negligible in the domain $[L+1,\infty]$, with $L = 2$. 
}

We observe that, if the potential does not depend on $\varepsilon$, $M \rightarrow 0$ as $\varepsilon \rightarrow 0$ and then the condition \eqref{BCt} reduces to a zero Neumann boundary condition. Therefore, in validating the model, we shall assume that $\varepsilon \to 0$ with $M$ fixed and we explain how to perform such a limit.

{  There are several ways to choose the values of the parameters $M$, $\varepsilon$, and of the function $U(\xi)$ (that describes the shape of the bubble interaction with the ions) in Eq.~\eqref{expr_M}. As we said before, since $M \to 0 $ such as $\varepsilon \to 0$, we may consider a potential $U = U^\varepsilon (\xi)$ which scales in such a way that $M$ is finite as $\varepsilon$ approaches $\vepst$, where $\vepst$ is a typical length of the attractive-repulsive core of the bubble and it is prohibitively small for a detailed calculation using Eq.~\eqref{equation_definition_flux}. In practice, we can assume that $M$ is non negligible when $\varepsilon\ll l$, where $l$ is a typical length scale of the physical setup. Here we show how to choose the parameters such that $M$ is constant when $\varepsilon \to \vepst$. The role of $\vepst$ will be evident in Sec.~\ref{section_saturation} where we describe the effect of saturation. 
	
	Once we define $U(\xi)$ in Eq.~\eqref{expr_U_LJ}, the resulting expression for $M$ from Eq.~\eqref{expr_M} is
	\begin{eqnarray}
		M &=& \varepsilon \mathcal{I}_L(\phi)
		\label{eq_choiseparam}
	\end{eqnarray}
	where 
	\begin{equation}
		\mathcal{I}_L(\phi) = \int_0^{L+1} \exp\left(-\phi\left( \xi^{-12} - 2\xi^{-6} \right) \right) d\xi.
		\label{eq_Iphi_mult}
	\end{equation}
	and it is represented in Fig.~\ref{fig:Iphi_mult}.
	Notice that, since the integrand in \eqref{eq_Iphi_mult} is exponentially large in $\phi$ over a finite interval, then 
	we can choose $(\varepsilon,\phi)$ such that $M$ is constant by solving the non linear equation for $\phi$, $\mathcal{I}(\phi) = M/\varepsilon$.
}
\begin{table}[]
	\begin{center}  
		\begin{tabular}{|c|c|c|c|c|}
			\hline
			$L$ & 2 & 4 & 6 & 8 \\ \hline
			ratio & $8\times 10^{-9}$ & $10^{-9}$ & $2\times 10^{-10}$ & $8\times 10^{-11}$ \\ \hline
		\end{tabular}
	\end{center}
		\caption{\textit{In this table we define {\rm ratio} := $ {\int_{L+1}^\infty U(\xi) d\xi}/{ \int_{0}^\infty U(\xi) d\xi }$, for different values of $L$, and $\phi = 6$. Here we show that we obtain a good approximation of the integral of the potential $U(\xi)$ with a few units of $L$.}}
		\label{table_L}
	\end{table}
	
	\subsection{  Extension to more dimensions}
	\begin{figure}[H]
		\centering
			\includegraphics[width=0.5\textwidth]{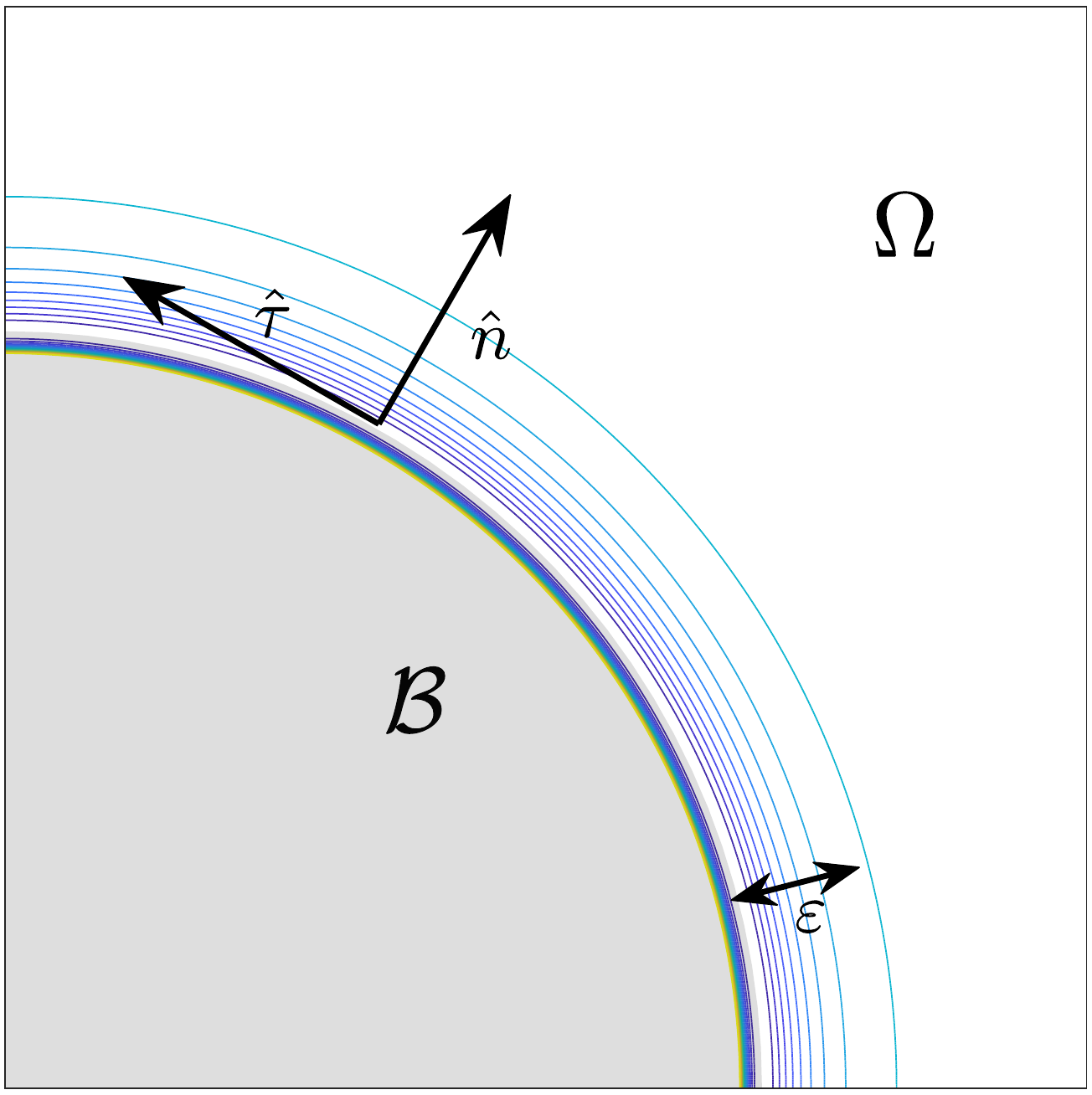}
			\caption{\textit{Scheme of the boundary condition in 2D where 
					$\hat{n}$ is the outgoing normal vector to the domain $\Omega$, {  represented by the white region outside the light gray bubble, and $\varepsilon$ is the thickness parameter of the bubble that may be a function of space. The contour lines represent the iso-potential lines. The white region near the bubble surface represents the local minima of the potential well, and represents the boundary $\Gamma_B$ of the bubble. The potential changes in the direction $\hat{n}$, while it is constant along the tangential direction $\hat{\tau}$.}}}
			\label{fig:bubble}
		\end{figure}
		
		{  In this section we provide a unified treatment of 2D and 3D case, making use of the projection operator on the line  (2D case) or the surface  (3D case), which delimits the bubble ${\mathcal B}$. 
			
			We assume that the bubble $\mathcal{B}$ is implicitly defined by a level set function $\Phi(\vx)$, $\vx\in\mathbb{R}^d$, $d=2,3$, and assume $\Phi(\vx)<0$ inside the bubble and $\Phi(\vx)>0$ in $\Omega$. For example $\Phi$ may be the signed distance from $\Sigma_B$ (see \cite{russo2000remark}), where $\Sigma_B$ denotes the boundary of the bubble $\mathcal{B}$. The unit normal $\hat{n}$ to $\Sigma_B$ is given by $\hat{n}(\vx) = \nabla\Phi(\vx)/|\nabla\Phi(\vx)|$, 
			$\forall \vx$ on $\Sigma_B$. Note that the $\hat{n}$ is naturally defined everywhere in $\Omega\cup\mathcal{B}$. We denote by $n_i, i=1,\ldots,d$ the Cartesian components of $\hat{n}$, and we shall denote by $h_{ij} = \delta_{ij}-n_in_j$ the projection operator on the tangent plane on $\Sigma_B$.  Finally, we denote by $\tpar_i \equiv h_{ij}\partial_j$ the gradient operator on the tangent plane to $\Sigma_B$. We adopt standard Einstein's convention of summation over repeated indices. Furthermore, we assume that the surface $\Sigma_B$ is smooth, and its curvature $\kappa$ times the range of the potential $\tilde{\veps}$ is much smaller than one.
			
			For points near the boundary, we decompose the flux vector $J_i,i=1,\ldots,d$ into a normal and tangential component, i.e.
			\begin{eqnarray}
				\label{eq_components_flux}
				J_i = \underbrace{n_i J_n}_{J^n_i} + \underbrace{h_{i j} J_j}_{J^\tau_i}
			\end{eqnarray}
			where $J_n = n_iJ_i$ and $J^\tau_i = h_{ij} J_j$.
			
			The divergence of the flux is given by:
			\begin{equation}
				\label{eq_flux_i}
				\nabla \cdot \vec{J} = \partial_i J_i
			\end{equation}
			If we substitute Eq.\eqref{eq_components_flux} in Eq.\eqref{eq_flux_i} we obtain
			\begin{equation}
				\nabla \cdot \vec{J} = \partial_i \left( n_i J_n \right) + \partial_i \left( h_{i j} J_j \right).
			\end{equation}
			Now the first term in the last equation can be written as
			\begin{align}
				\label{eq_n_i}
				\partial_i \left( n_i J_n \right) & = n_i \partial_i J_n + J_n \underbrace{\partial_i n_i}_{\chi_{ii}} &
			\end{align}
			where $\chi_{ij} = \partial_in_j$ denotes the second fundamental form of the surface  $\Sigma_B$, and its trace $\chi_{ii}$ gives the curvature of the surface $\Sigma_B$:
			\begin{eqnarray}
				\chi_{ii} = \kappa = \left\{
				\begin{array}{rl} \displaystyle
					1/R &  {\rm in}\, 2D \\
					\kappa_1 + \kappa_2 &  {\rm in} \, 3D, 
				\end{array}
				\right.
			\end{eqnarray}
			where $R$ is the local radius of curvature,  $\kappa_1$ and $\kappa_2$ denote the Gauss principal curvatures.
			
			Since the gradient of the potential is parallel to the normal, then one has
			\begin{equation}
				\label{eq:Jtau}
				J^\tau_i = -D\tpar_i c = -D h_{ij}\partial_j c
			\end{equation}
			At this point we substitute the quantities (\ref{eq_n_i}-\ref{eq:Jtau}) in Eq.\eqref{eq_flux_i} and the expression for the divergence of the flux becomes:
			\begin{equation}
				\nabla \cdot \vec{J} = n_i \partial_i J_n + J_n\chi_{ii} - D\Delta_\perp c
			\end{equation}
			where $\Delta_\perp c = \tpar_i\tpar_i c = h_{ij}\partial_i\partial_j c - \chi_{ii}n_j\partial_j c $ denotes the surface Laplacian of the concentration (see Appendix~\ref{A:LB}).
			
			The evolution equation of for the concentration near the bubble surface is then given by
			\begin{equation}
				\pad{c}{t} = -\pad{J_n}{n} - J_n \chi_{ii} + D \Delta_\perp c.
			\end{equation}
			
			Following the procedure adopted in 1D, let $\vec{x}$ be a point in $\Sigma_B$ and integrate along the normal direction:
			\begin{align} \nonumber
				\frac{\partial}{\partial t}\int_{-\varepsilon}^{\varepsilon L} c (\vec{x}+r\hat{n})\,dr &= -J_n(\vec{x}+L\varepsilon\hat{n}) 
				- \int_{-\varepsilon}^{\varepsilon L} \chi_{ii}  J_n (\vec{x}+r\hat{n}) \, dr  \\  
				& + D \Delta_\perp \int_{-\veps}^{\veps L}  c (\vec{x}+r\hat{n})\,dr \nonumber
			\end{align}
			{ 
				where we dropped the term $J_n(\vec{x}-\veps\hat n)$ which vanishes because of the repulsive core. 
				
				Notice that the term 
				$\int_{-\varepsilon}^{\varepsilon L} \chi_{ii}  J_n (\vec{x}+r\hat{n}) \, dr  = (1+L)\veps\langle \chi_{ii} J_n \rangle$ can be neglected since we assume that the range of the potential is much smaller than the radius of curvature (2D) or to the inverse of the mean curvature (3D).
			}

			Now we consider a change of variable from $\vec{x}$ to $\xi$ in the integration intervals, and obtain 
			\begin{align} \nonumber
				\veps\frac{\partial}{\partial t}\int_0^{L+1} c (\vec{x}+\veps\xi\hat{n})\,dr &= -J_n(\vec{x}+L\varepsilon\hat{n}) + \veps D \Delta_\perp \int_{0}^{L+1}  c (\vec{x}+\veps \xi \hat{n})\,d\xi\nonumber
			\end{align}
			
			Let us consider {  the normal components of the flux $J_n$}:
			\begin{equation}
				\label{eq_J_n}
				J_n = -D\left(\frac{\partial  c }{\partial {n}}+\frac{1}{k_B T}\frac{\partial V}{\partial {n}} c\right),
			\end{equation}
			since $\displaystyle V(x)/k_B T=U\left(\xi\right)$, {  it means that we have the following relation between the two expressions of the potentials }
			\[
			\displaystyle \frac{\partial (V/k_B T)}{\partial {n}} = \frac{1}{\varepsilon}U'(\xi)
			\] 
			{   We notice, as in the 1D case, that} there is also a dependence on $\xi $ of the solution inside the layer, such that $ c = c(\vec{x} + \varepsilon \xi \hat{n},t) = \tilde{c}(\xi,t)$, thus we can put $1/\varepsilon$ as common factor:
			\begin{equation}
				J_n = -D\frac{1}{\varepsilon}\left(\frac{\partial  c }{\partial \xi}+ \frac{\partial U}{\partial \xi}c\right)
			\end{equation}
			where we omitted the $\tilde{\cdot}$. Following the same argument we used in 1D to derive Eq.~\eqref{exprho}, we perform a formal expansion in $\varepsilon$ of the solution. To the lowest order in $\varepsilon$ {  for the normal component of the flux we have $\displaystyle {J_n^{(0)}} = 0$, } which implies 
			\begin{equation}
				\label{bcbubble}
				\frac{\partial  c^{(0)}}{\partial \xi} + \frac{\partial U}{\partial \xi} c^{(0)} = 0.
			\end{equation}
			Therefore we obtain {  the expression for the concentration to order zero in $\varepsilon$:}
			\begin{equation}
				c^{(0)}(\vec{x}+\varepsilon \xi\hat{n})= c^{(0)}(\vec{x}+\varepsilon L\hat{n})\exp\left(-U(\xi)\right).
			\end{equation}
			Using this expression, we can compute the line density (2D) or the surface density (3D) of entrapped ions:
			\begin{equation}
				\mathcal{C}(\vec{x}) = \varepsilon\int_{0}^{L+1} c (\vec{x}+\varepsilon \xi\hat{n})d\xi \simeq M c (\vec{x})
			\end{equation}
			with $\displaystyle M= \varepsilon \int_{0}^{L+1}\exp\left(-U(\xi)\right)d\xi$ and where we assumed $U(L+1)=U'(L+1)=0$.
			
			{ At the end, in 2D (3D),} the system and the new boundary condition for the concentration becomes:
			\begin{eqnarray}
				\label{system_multiscale}
				\displaystyle \frac{\partial c}{\partial t} &=& D \Delta c \quad \rm{in} \> \Omega \\ \nonumber
				\displaystyle \frac{\partial c}{\partial n} &= &0 \quad \rm{ on } \, \partial \Omega\backslash \Sigma_B \\ \label{bc_multiscale}
				\displaystyle M\frac{\partial  c }{\partial t} &=& MD\Delta_\perp  c +D\frac{\partial  c }{\partial n} \quad \rm{ on }\, \Sigma_B
			\end{eqnarray}
			because $J_n = -D\partial  c /\partial n$ just out of the range of the potential.
			
			Notice that in 2D the Laplace-Beltrami operator reduces to the second derivative with respect to the arclength of the boundary:\footnote{Here by $\partial/\partial \tau$ 
				we denote the derivative on $\Gamma_B$, i.e.\ the derivative along the arclength that parametrizes the curve, likewise
				$\partial^2/\partial \tau^2$ denotes the second derivative along $\Gamma_B$, not the second derivative along the tangent direction. See Fig.\ \ref{fig:bubble} and Appendix~\ref{A:LB}. }
			\begin{equation}
				\Delta_\perp c = \frac{\partial^2 c}{\partial \tau^2}.
				\label{eq:BL2D}
			\end{equation}
		}
		
\section{Validation of the model and results}
\label{section_validation_multiscale}
In this section we validate the model in 1D and 2D, considering different geometries. {  In all numerical tests, when dealing with the full model, we choose a space discretization that fully resolves the $\varepsilon$-scale}.

We start considering the effect of the LJ potential in 1D and the agreement with the multiscale model. To validate the model we compare the mass entrapped in the hole {  for the full model,} with the concentration calculated in the {  boundary} and compare the profiles of the fully resolved solution and of the solution of the multiscale model in the bulk. 

Secondly, we consider a 2D system in a squared domain with a potential that is constant in the $y$-direction and concentrated on the left edge of the domain.

Finally we consider a more realistic domain in 2D, with a circular hole describing the attractive bubble surface at the center of a squared domain. 

\subsection{Validation 1D}
The domain in 1D is $\Omega^\varepsilon= [- \varepsilon,1] =\Omega^\varepsilon_b \cup \Omega^\varepsilon_f ,\,\Omega^\varepsilon_b = [- \varepsilon,L\varepsilon],\, \Omega^\varepsilon_f = [L\varepsilon,1]$. {  In all our calculations we used $L=2$, as we explain in Section \ref{section_1Dmodel}.} The computational domain $\Omega^\varepsilon_h$ is then discretized by a uniform Cartesian mesh with spatial step $ h $. The concentration $c_{i}$ is defined at the center of the cells with $x_i=-\varepsilon + (i-1/2) h \in\{1,\dots,N_x\},  h  N_x = 1 +  \varepsilon$. We choose a cell centered discretization in order to guarantee the exact conservation of the  
{total volume of the solute $v=\sum_i c_i h$} deriving by the zero boundary condition for the flux. The scheme is second order accurate and it is stable even in presence of a drift term, provided the so called \emph{mesh P\'{e}clet number} is smaller than 2 \cite{Wesseling2023600}, defined as follows:
\begin{equation}
	\label{eq_pec}
	{\rm pec} \equiv \max_x\left|{\partial_x U}\right| h < {2}.
\end{equation} 

For the multiscale model the domain is $\Omega^0 = [0,1]$, and again we discretize the computational domain $\Omega^0_h$ as follows: $x_i=(i-1/2) h \in\{0,\dots,N\},  h  N = 1$. Internal equations are discretized by finite differences by the classical three points approximation of the second derivative, while the Eq.~\eqref{BCt} is discretized to second order as 
\[ 
M\frac{\partial }{\partial t} \frac{c_0 + c_1}{2} = D \frac{c_1 - c_0}{h}. 
\]
In this case we do not have restriction defined in Eq.~\eqref{eq_pec} on the space step because $U=0$ in $\Omega^0$.

If we denote with $c^\varepsilon$ and $c^0$ the solutions of the full and multiscale models, respectively, the problems~\eqref{pde1d}-\eqref{flux1d} and \eqref{reduced1d}-\eqref{expr_M} can be summarized as 
\begin{equation}\label{compactQ}
	\frac{\partial  c^\varepsilon }{\partial t} = Q^\varepsilon ( c^\varepsilon ),
\end{equation}
where we distinguish between the cases $\varepsilon > 0$ and $\varepsilon = 0$:
\begin{eqnarray}
	\label{eq_Q_eps}
	Q^\varepsilon( c^\varepsilon ) =	\left\{
	\begin{matrix} 
		D \Delta  c^\varepsilon  & {\rm in } \, \Omega^\varepsilon_f\\
		\displaystyle D \left(\Delta  c^\varepsilon + \frac{1}{k_B T}\nabla \cdot (c^\varepsilon V')\right) & \rm{ in } \,\Omega^\varepsilon_b 
	\end{matrix}
	\right. \qquad {\rm if} \> \varepsilon>0\\
	\label{eq_Q_0}
	Q^0( c^0 ) =	\left\{
	\begin{matrix} 
		D \Delta  c^0  & {\rm in } \, \Omega^0, \\
		\displaystyle DM^{-1} \frac{\partial  c^0 }{\partial x} \quad & {\rm on }\, x=0 
	\end{matrix}
	\right. \qquad {\rm if} \> \varepsilon=0   
\end{eqnarray}
Eq.~\eqref{compactQ} is discretized in time by using the \textit{Crank-Nicolson} method, which is second order accurate:
\begin{eqnarray}
	\frac{ c  ^{\varepsilon,n+1}- c ^{\varepsilon,n}}{\Delta t} &=& \frac{1}{2}Q^\varepsilon( c ^{\varepsilon,n}) + \frac{1}{2}Q^\varepsilon( c ^{\varepsilon,n+1})
	\label{CNdisc}
\end{eqnarray}
where $\Delta t$ is the time step, and we solve the system for $c^{\varepsilon,n+1}$. Eq.~\eqref{compactQ} is supplemented by the following initial condition:
\begin{equation}
	c(t=0,x)=c_{\rm in}(x)=\frac{v_0}{\sqrt{2\sigma^2\pi}}\exp\left(-\left(x-x_m\right)^2/2\sigma^2\right), \qquad 
	\sigma \in \mathbb{R}, \quad x_m
	\in(0,1),
	\label{ic_1d}
\end{equation}
and $v_0$ denotes the total volume. 

\begin{figure}[H]
	\begin{minipage}
		{.44\textwidth}
		\centering
		\includegraphics[width=\textwidth]{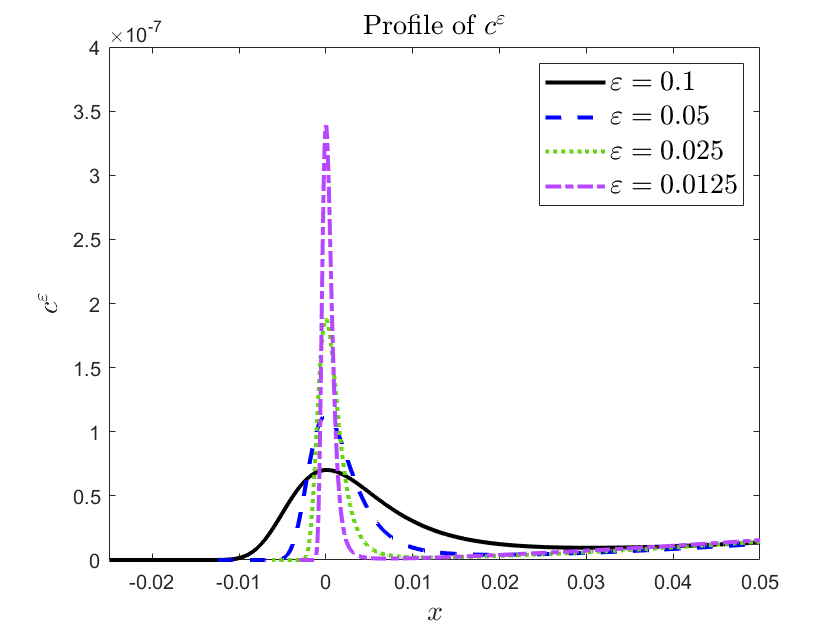}
	\end{minipage}
	\begin{minipage}
		{.44\textwidth}
		\includegraphics[width=\textwidth]{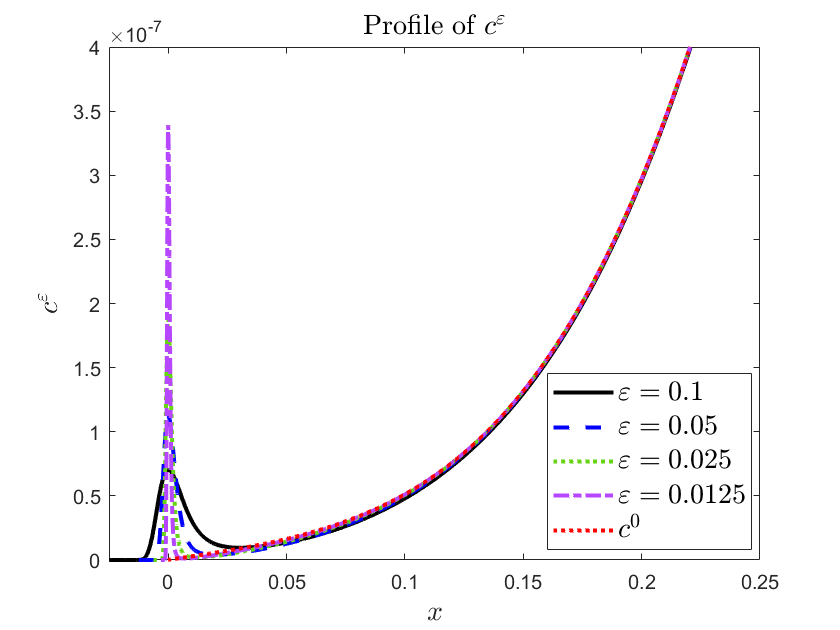}
	\end{minipage}
	\caption{\textit{Here we see two aspects of the same plot. On the left we see the effect of the attractive potential that becomes stronger decreasing the value of $\varepsilon$. We consider a zoom-in in $x$-direction to focus on the length of the entrapped mass. In the plot on the right we have a zoom-in in $y$-direction with an additional case. The dashed red line represents the solution of the multiscale model, $c^0$. In this panels we see the agreement of the two models improves with $\varepsilon \to 0$. Parameters of the tests: $\Delta x = 10^{-4},\,t_{\rm final} = 0.05,\,M=3$  and for the initial condition defined in Eq.~\eqref{ic_1d} we choose $\sigma = 0.1$, $v_0 = 10^{-6}$, $x_m = 0.5$. Parameters of the tests: $\Delta x = 10^{-4},\,t_{\rm final} = 0.05,\,M=3$ and for the initial condition defined in Eq.~\eqref{ic_1d} we choose $\sigma = 0.1$, $v_0 = 10^{-6}$, $x_m = 0.5$.}}
	\label{1D_full_complete}
\end{figure}

The method that we consider takes into account two parameters, the range of the potential $\varepsilon$ and the discretization step $h=\Delta x$ (we chose $\Delta t \propto \Delta x$). In this section first we see the effect of the LJ potential, described in Eq.~\eqref{expr_U_LJ}, on the concentration (see Fig.~\ref{1D_full_complete}, left). As $\varepsilon$ decreases, the agreement between the two models gets better and better, as we see in Fig.~\ref{1D_full_complete} (right).

In order to perform a quantitative comparison between the two models, we calculate the difference in the amount of solute near the surface, which we denote by $e_S$, and the norm of the concentration difference in the bulk, denoted by $e_B$ (scaled by the total amount of the solute $v_0$), as functions of time:
\begin{eqnarray}
	\label{eq_errors_time}
	\displaystyle {\rm e_S} = \frac{\left|\int_{-\varepsilon}^{\varepsilon L}c^\varepsilon(x,t)d x - M c^0(x=0,t) \right|}{v_0} , \quad {\rm e_B }= \frac{\int_{\varepsilon L}^{1} \left|c^\varepsilon(x,t) - c^0(x,t)\right|d x }{v_0}  \\ \label{eq_error}
	{\rm error} = \frac{|\int_{-\varepsilon}^{L\varepsilon} c^\varepsilon(x)d x - M c^0(x=0,t_{\rm final})|}{|M c^0(x=0,t_{\rm final})|} 
\end{eqnarray}

\begin{figure}[H]
	\centering
	\begin{minipage}
		{.48\textwidth}
		\centering
		\includegraphics[width=\textwidth]{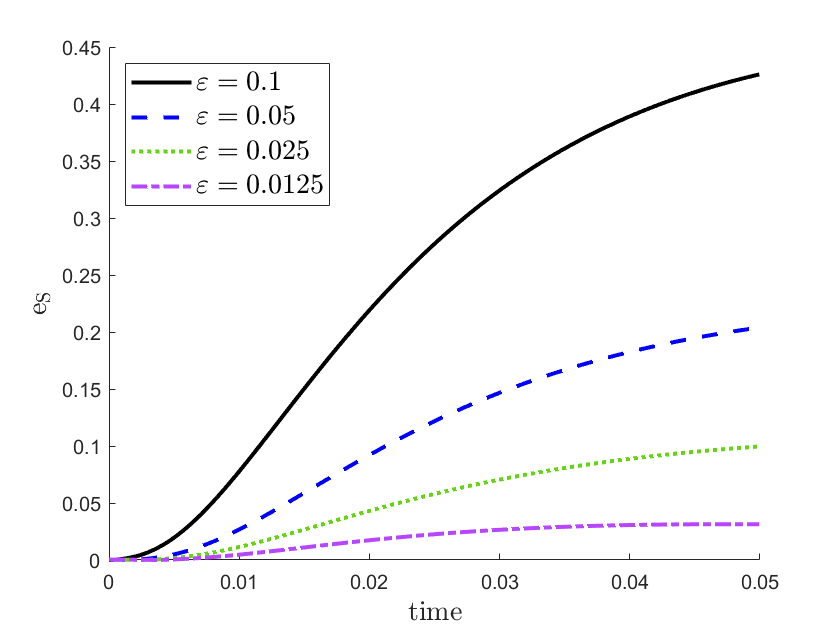}
	\end{minipage}
	\begin{minipage}
		{.48\textwidth}
		\includegraphics[width=\textwidth]{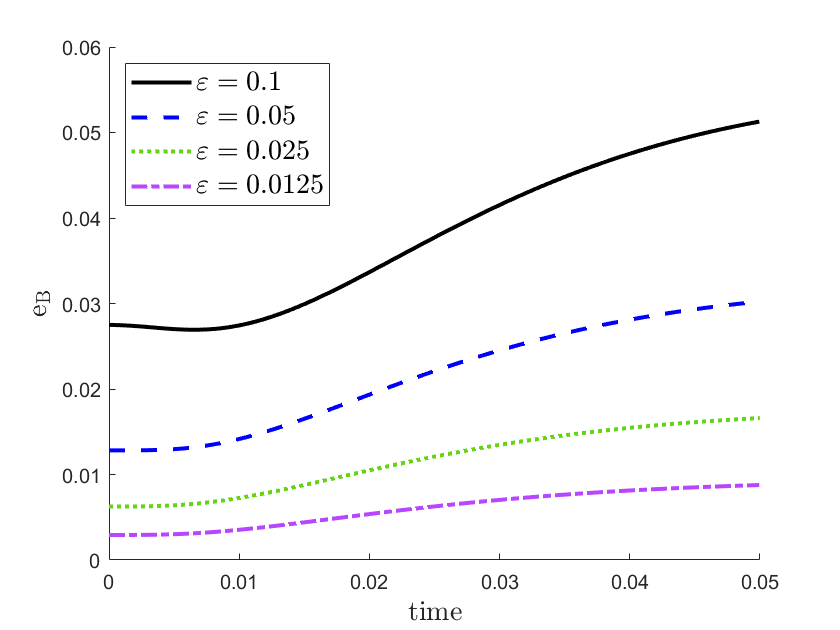}
	\end{minipage}
	\caption{\textit{Errors defined in Eq.~\eqref{eq_errors_time} for the entrapped mass near the surface (left panel) and in the bulk (right panel). We see the errors decrease with $\varepsilon$ as expected. Parameters of the simulations: $\Delta x = 1.82\times 10^{-4},\,t_{final} = 0.05$ and for the initial condition defined in Eq.~\eqref{ic_1d}: $\sigma = 0.2$, $v_0 = 10^{-6}$, $x_m = 0.5$.}}
	\label{errors_mass_time}
\end{figure}

The results are reported in 
Fig.~\ref{errors_mass_time}. As we can see, the discrepancy between the two approaches is approximately proportional to the range of the potential $\varepsilon$.

\begin{figure}[H]
	\centering
	\begin{minipage}
		{.48\textwidth}
		\centering
		\includegraphics[width=\textwidth]{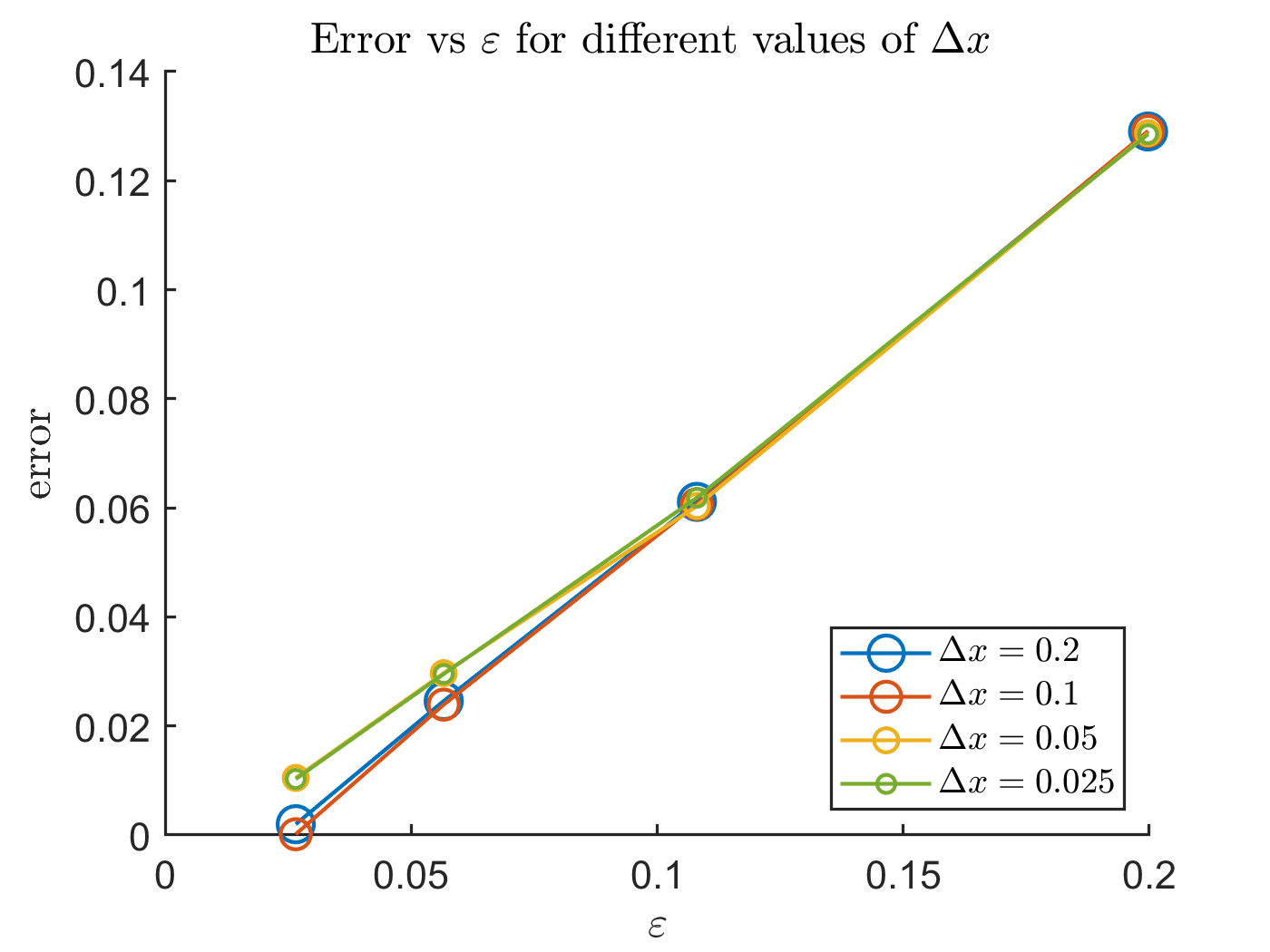}
	\end{minipage}
	\begin{minipage}
		{.48\textwidth}
		\includegraphics[width=\textwidth]{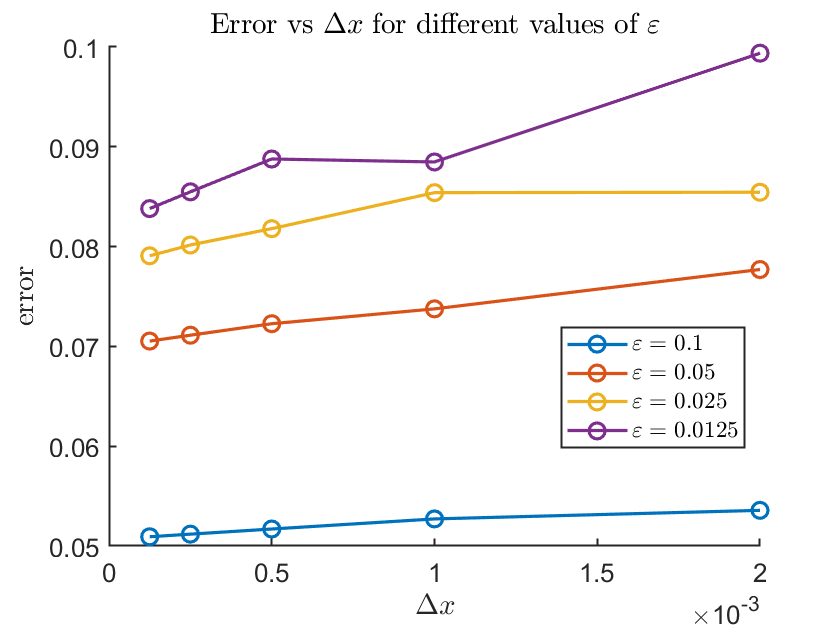}
	\end{minipage}
	\caption{\textit{In this figure we show how the error, defined in Eq.~\eqref{eq_error}, does not depend on the spatial discretization, once the $\varepsilon$-scale is resolved. It also shows that error increase linearly with $\varepsilon$. This is why we have the same behaviour for all values of $\Delta x$ (left panel) and the errors are constants as functions of $\Delta x$ (right panel).}}
	\label{error_1D_deltax}
\end{figure}

Finally in 1D we  {show} that the error between models does not depend on the space discretization. In Fig.~\ref{error_1D_deltax} we see the error is  {almost} independent of the space step $\Delta x$,  {especially when $\Delta x \ll \veps$.}

\begin{figure}[H]
	\begin{minipage}
		{.44\textwidth}
		\includegraphics[width=\textwidth]{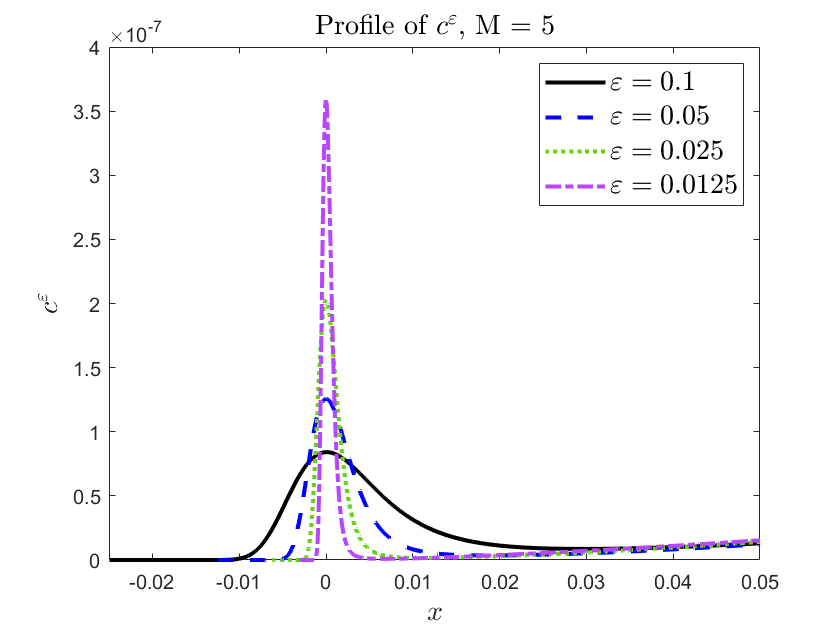}
	\end{minipage}
	\begin{minipage}
		{.44\textwidth}
		\includegraphics[width=\textwidth]{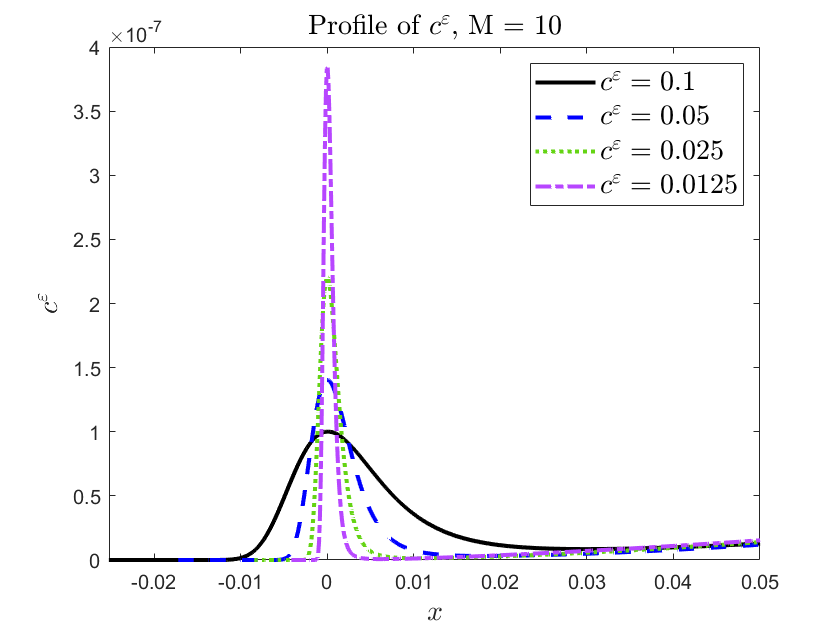}
	\end{minipage}
	\caption{\textit{{  Here we see the effect of the potential with two different values of $M$, on the left $M = 5$ and on the right $M=10$. We see that for a higer value of $M$, (but same values of $\varepsilon$) the peak of the solution is higher. Parameters of the tests: $\Delta x = 10^{-4},\,t_{\rm final} = 0.05,$ and for the initial condition defined in Eq.~\eqref{ic_1d} we choose $\sigma = 0.1$, $v_0 = 10^{-6}$, $x_m = 0.5$.} }}
	\label{1D_full_complete_2}
\end{figure}	


\subsection{Validation in 2D}
\label{sec:2D_easy}
In this section we define a very simple domain in 2D. It is a direct extension of the 1D case, with the potential defined as a function $V:=V(x,y)$ but it is constant in $y$-direction. In this case the domain is defined as $\Omega^\varepsilon = [-\varepsilon,1]\times[0,1]$ and the potential is negligible for $x>L\varepsilon $ with $L=2$ as we do in 1D. 

The computational domain $\Omega^\varepsilon_h$ is then discretized by a uniform Cartesian mesh with spatial steps $\Delta x$ and $\Delta y$ in $x$ and $y$ directions, respectively. The concentration $c_{ij}$ is defined at the center of the cells with $x_i=- \varepsilon + (i-1/2)\Delta x\in\{1,\dots,N_x\}, \Delta x N_x = 1 +  \varepsilon$ and $y_j=(j-1/2)\Delta y\in\{1,\dots,N_y\}, \Delta y N_y = 1$. We choose again a cell centered discretization, and  {we make sure that the stability condition imposed by the \emph{mesh P\'{e}clet number} always holds.}

For the multiscale model the domain is $\Omega^0 = [0,1]\times[0,1]$, and again we discretize the computational domain $\Omega^0_h$ as follows: $x_i=(i-1/2)\Delta x\in\{0,\dots,N\}, \Delta x N = 1$ and $y_i=(i-1/2)\Delta x\in\{1,\dots,N\}, \Delta y N = 1$, and $h = \Delta x = \Delta y$.

\begin{figure}[H]
	\centering	
	\begin{minipage}
		{.48\textwidth}
		\centering
		\includegraphics[width=\textwidth]{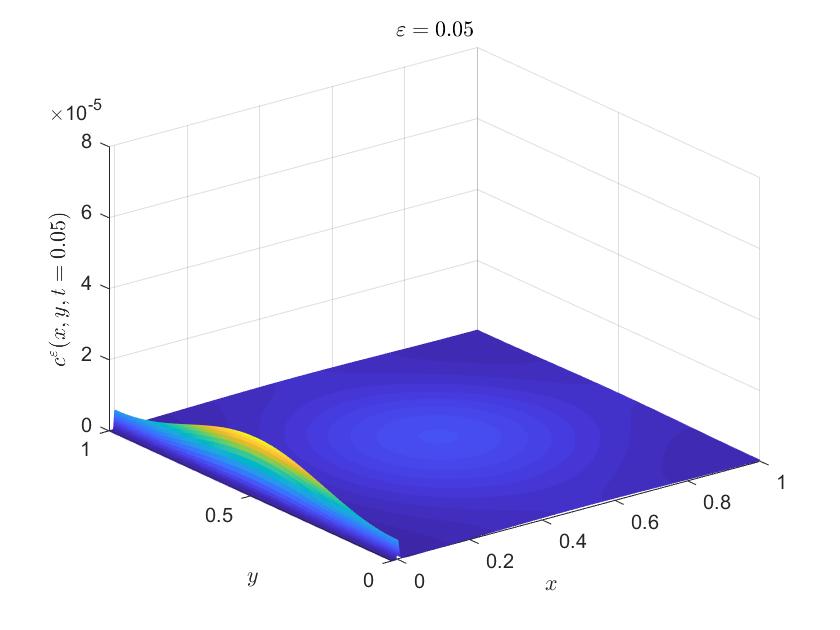}
	\end{minipage}
	\begin{minipage}
		{.48\textwidth}
		\centering
		\includegraphics[width=\textwidth]{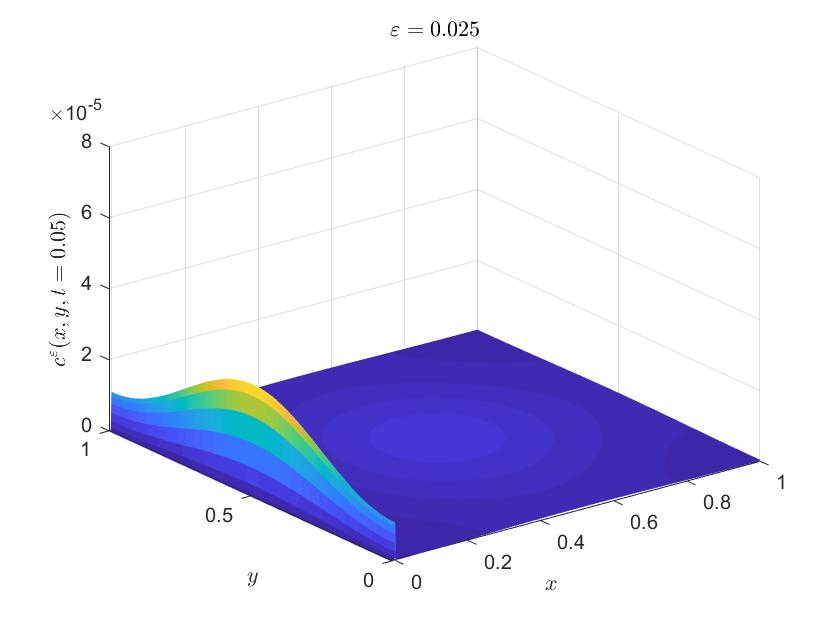}
	\end{minipage}
	\begin{minipage}
		{.48\textwidth}
		\centering
		\includegraphics[width=\textwidth]{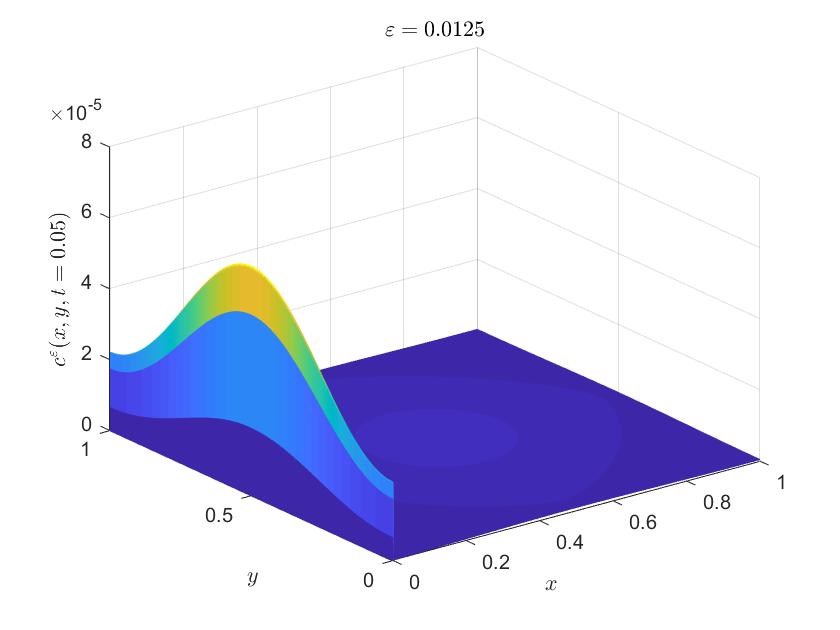}
	\end{minipage}
	\begin{minipage}
		{.48\textwidth}
		\centering
		\includegraphics[width=\textwidth]{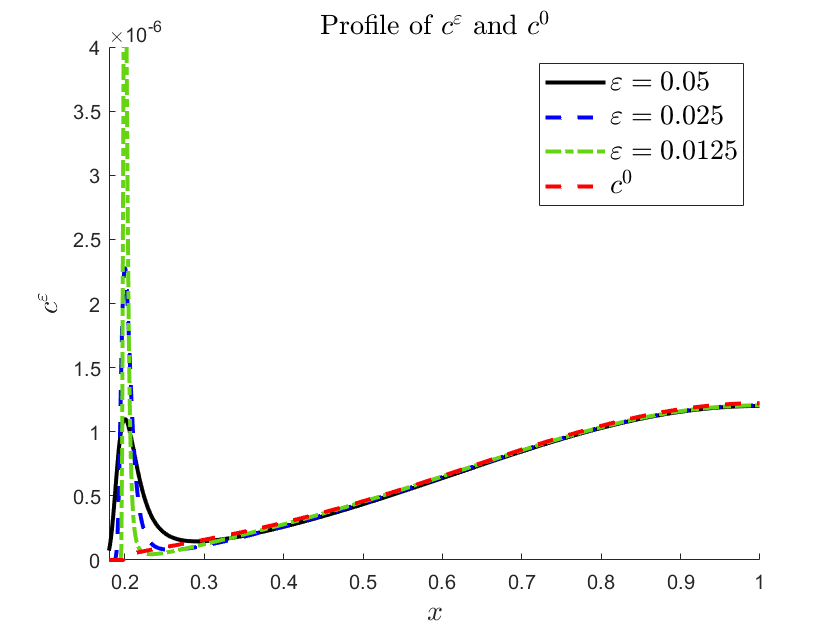}
	\end{minipage}
	\caption{\textit{We see the effect of the attractive potential that becomes stronger and stronger decreasing the value of $\varepsilon$ in panels (a,b,c). In panel (d) we have a 1D plot, that is a section of the previous one for $y=0.5$ with an additional line. The dashed red line represents the solution of the multiscale model, $c^0$. The agreement of the two models improves with $\varepsilon \to 0$. Grid resolution of the simulations: $\Delta x = \Delta t = 3.33\times 10^{-4}$, $t_{final} = 0.05s$ and the initial condition is defined in Eq.~\eqref{ic_1d} with $x_m = y_m = 0.5, \sigma = 0.1, v_0 = 10^{-6}$. }}
	\label{figure_effect_epsilon_2D}
\end{figure}

{About the discretization in time, in both cases, full and multiscale model, we choose the Alternate Direction Implicit (ADI) method \cite{Wesseling2023600}, which is basically as accurate as Cranck-Nicholson, but it is much mode efficient, since it requires to solve only tri-diagonal systems at each time step.
}
In Fig.~\ref{figure_effect_epsilon_2D} panels (a,b,c) we see the effect of the potential $V$ as $\varepsilon\to 0$ and in panel (d) we see a one dimensional plot as a summary of the three previous cases ( {black, blue and green} lines) and the solution of the multiscale model (dashed red line).


\subsection{Validation 2D in presence of a circular bubble}
\label{sec:validation2D}
Here we consider a more realistic 2D domain, in which the bubble is described as a circular disk centered at the origin, for the multiscale model, while for the full problem the potential is defined as a function  $V(x,y):=\mathcal{U}\left({(r-R)}/{\varepsilon}\right)$ where $\mathcal{U}$ is a function $\mathcal{U}\colon \mathbb{R}\to\mathbb{R}$ such that
\[
\displaystyle \mathcal{U}(\xi)=\left\{
\begin{array}{l}
	+\infty \quad \xi<-1\\
	U(\xi) \quad \xi \in [-1,L] \\
	0 \quad \quad \xi>L
\end{array}
\right.
\]
with 
\[
\displaystyle \xi = \frac{r - R}{\varepsilon} = \frac{\sqrt{x^2 + y^2} - R}{\varepsilon}
\]
where $R$ is the radius of the bubble and $r$ is the distance from the origin,
and $U(\xi)$ is given by Eq.~\eqref{expr_U_LJ}.
The equation for  {$c=c^0$ in} the domain $\Omega^0$ is the following
\begin{equation}\label{pde2d}
	\frac{\partial   c }{\displaystyle \partial t} = \nabla \cdot \left( D \nabla  c  \right) \text{ in } \Omega^0.
\end{equation}
The fluid is contained in a box $\mathcal{S} = (-a,a)^2 \subseteq \mathbb{R}^2$, with $a>0$. Wall boundary conditions are prescribed on the boundary $\Gamma_\mathcal{S} = \partial \mathcal{S}$ (zero Neumann conditions for the solution):
\begin{equation}\label{bcwall}
	\frac{\partial  c }{\partial {n}} = 0 \quad \text{ on } \Gamma_\mathcal{S}.
\end{equation}
\begin{figure}[H]
	\centering
	\hfill
	\begin{minipage}
		{.45\textwidth}
		\centering
		\includegraphics[width=0.8\textwidth]{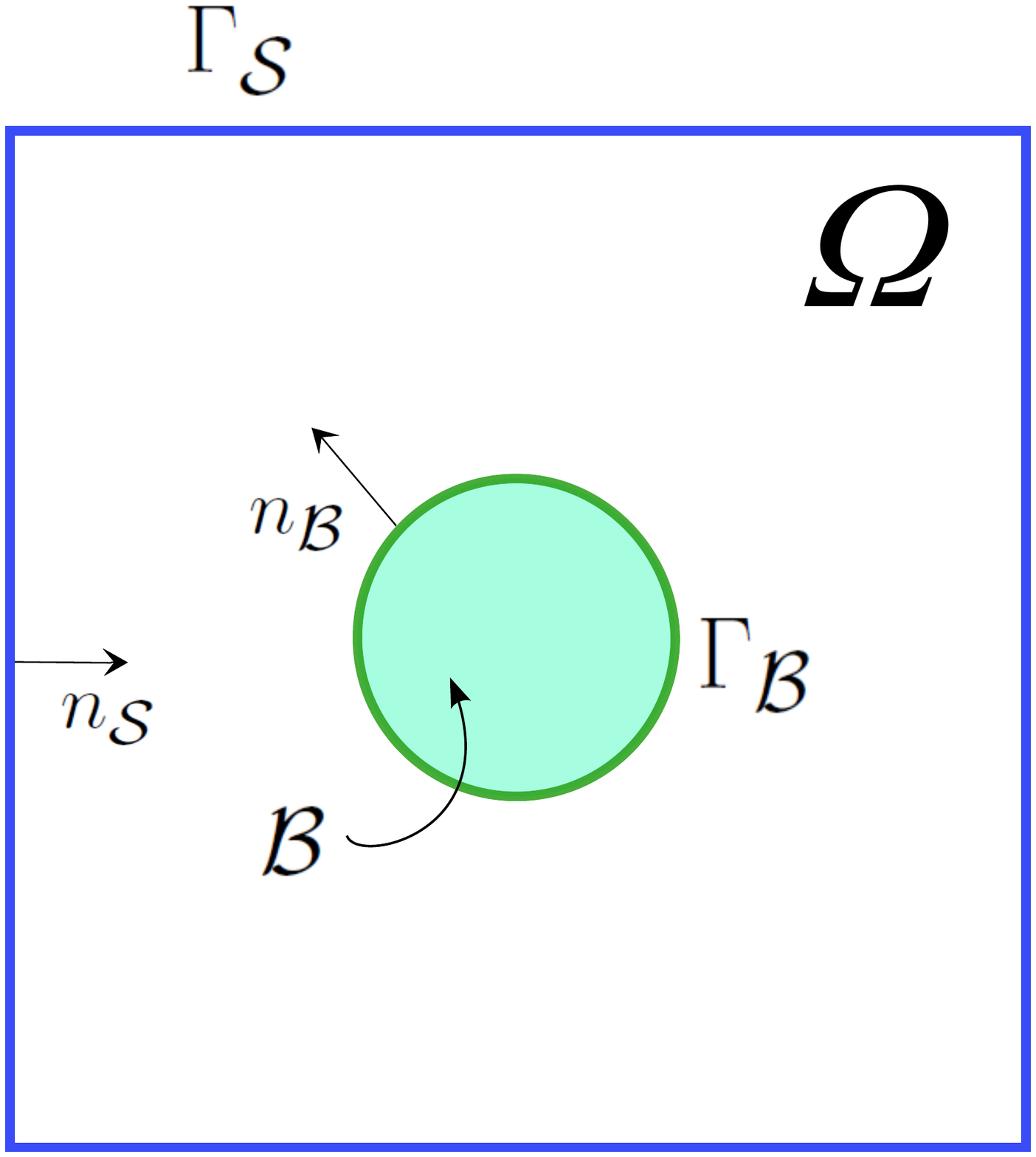}
	\end{minipage}\hfill
	\begin{minipage}
		{.45\textwidth}
		\centering
		\includegraphics[width=0.55\textwidth]{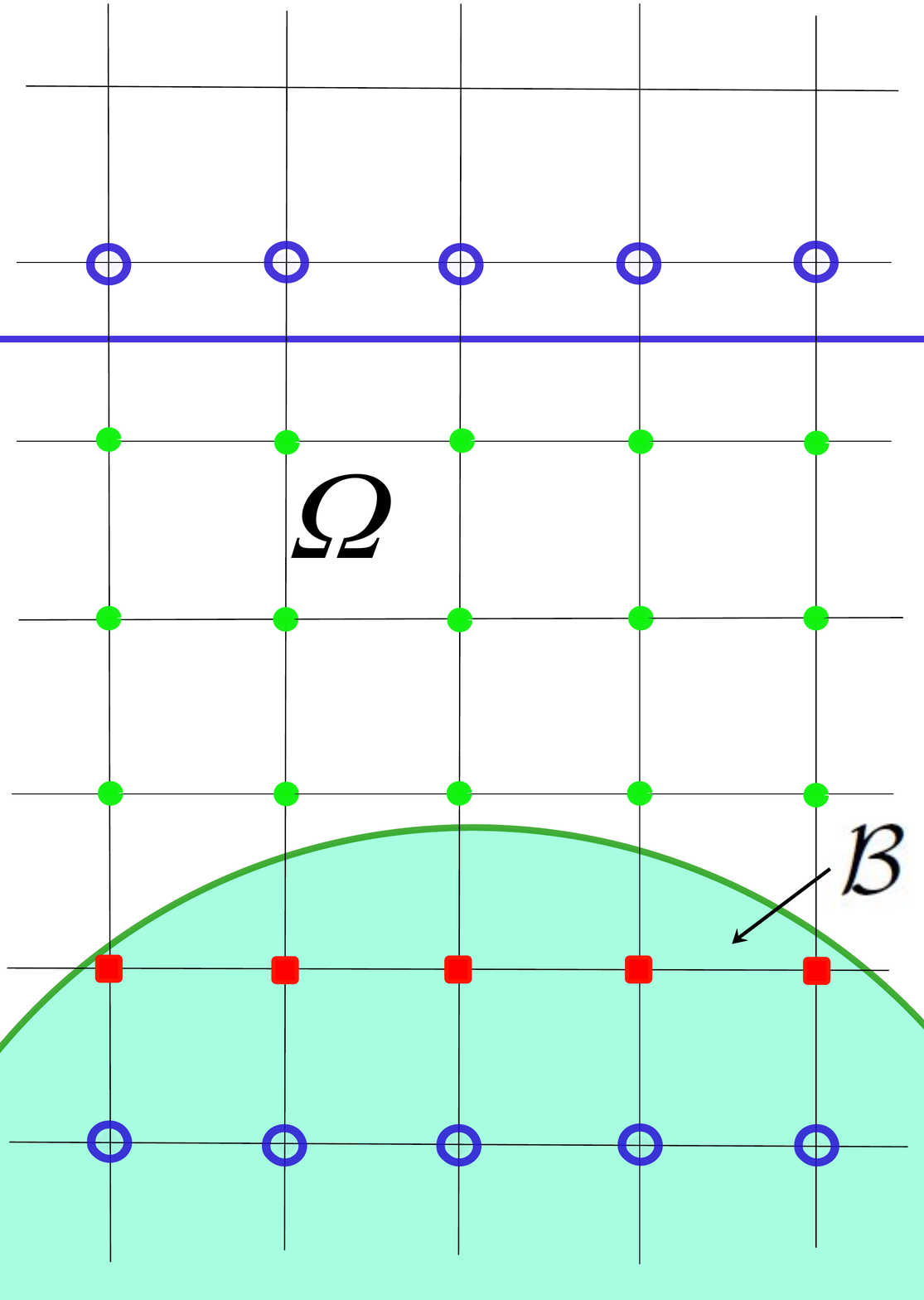}
	\end{minipage}
	\hspace*{\fill}
	\caption{\textit{Representation of the domain on the left and classification of inside grid points (green points), ghost points (red squares) and inactive points (blue circles) on the right.}}
	\label{classification_points}
\end{figure}
About the discretization in space for the full problem, it is analogue to the previous case, described in \ref{sec:2D_easy} and the expression for the initial condition is the following:
\begin{equation}
	c(t=0,x,y)=c_{\rm in}(x,y)=\frac{2v_0}{\sqrt{2\sigma^2\pi}}\exp\left(\left(-\left(x-x_m\right)^2+\left(y-y_m\right)^2\right)/2\sigma^2\right), \quad 
	\sigma \in \mathbb{R}
	\label{ic_2d}
\end{equation}

In the presence of a bubble, the fluid domain is represented by $\Omega^0 = \mathcal{S} \backslash \mathcal{B}$, where $\mathcal{B}$ is the region occupied by the bubble and represented by a circle centred in the origin and with radius $R$ such that $0 < R < a$. We choose $(x_m,y_m)\in \Omega^0$.

The computational domain $\mathcal{S}$ is discretized through a uniform Cartesian mesh with spatial step $h = 2a/N = \Delta x = \Delta y$, where $N^2$ is the number of cells.
We use a cell-centered discretization, therefore the set of grid points is $\mathcal{S}_h = \{(x_i,y_j)=(-a - h/2 + ih,-a - h/2 + jh), (i,j) \in \{1,\cdots,N\}^2 \}$. Within the set of grid points we define the set of internal points $\Omega^0_h = \mathcal{S}_h \cap \Omega^0$, the set of bubble points $\mathcal{B}_h =\mathcal{S}_h \cap \mathcal{B}$ and the set of ghost points $\mathcal{G}_h$:
\begin{equation}
	(x_i,y_j) \in \mathcal{G}_h \iff (x_i,y_j) \in \mathcal{B}_h \text{ and } \{(x_i \pm h,y_j),(x_i,y_j\pm h) \} \cap \Omega^0_h \neq \emptyset.
\end{equation}
The remaining grid points that are neither inside nor ghost points are called inactive points. See Fig.~\ref{classification_points} 
(right panel) for a classification of inside, ghost, and inactive points.

For an accurate description of the discretization of the problem we refer to \cite{astuto2023finite}, where we present a ghost-point method to solve drift-diffusion equations in domains with circular holes, in 2D and 3D, with boundary conditions defined in Eq.~\eqref{expr_M}, even in presence of a moving bubble. A suitable geometric multigrid approach is adopted to solve the problem, based on a proper relaxation of the boundary conditions on the ghost points, with a relaxation parameter chosen in order to achieve the convergence of the iterative method \cite{COCO2013464}. Here we briefly describe the space and time discretization.

The problem \eqref{CNdisc} is discretized in space, leading to a linear system
\begin{equation}\label{CNdiscspace}
	\left(I_h - \frac{k}{2}Q_h\right)  c _h^{n+1} = \left(I_h + \frac{k}{2}Q_h\right)  c _h^n
\end{equation}
to be solved at each time step, where $I_h$ and $Q_h$ are the $(N_I+N_G) \times (N_I+N_G)$ matrices representing the discretization of the identity operator $I$ and the operator $Q$ in Eqs.~(\ref{eq_Q_eps}-\ref{eq_Q_0}).

{While for the full model we use an ADI discretization in time, for the multiscale model we use a Crank-Nicolson method. The reason of this choice comes from the shape of $\Gamma_B$, that makes the discretization in time depending on the coordinates, losing the benefits of the ADI scheme.}
\begin{figure}[H]
	\centering
	\hfill
	\begin{minipage}
		{.45\textwidth}
		\centering
		\includegraphics[width=0.85\textwidth]{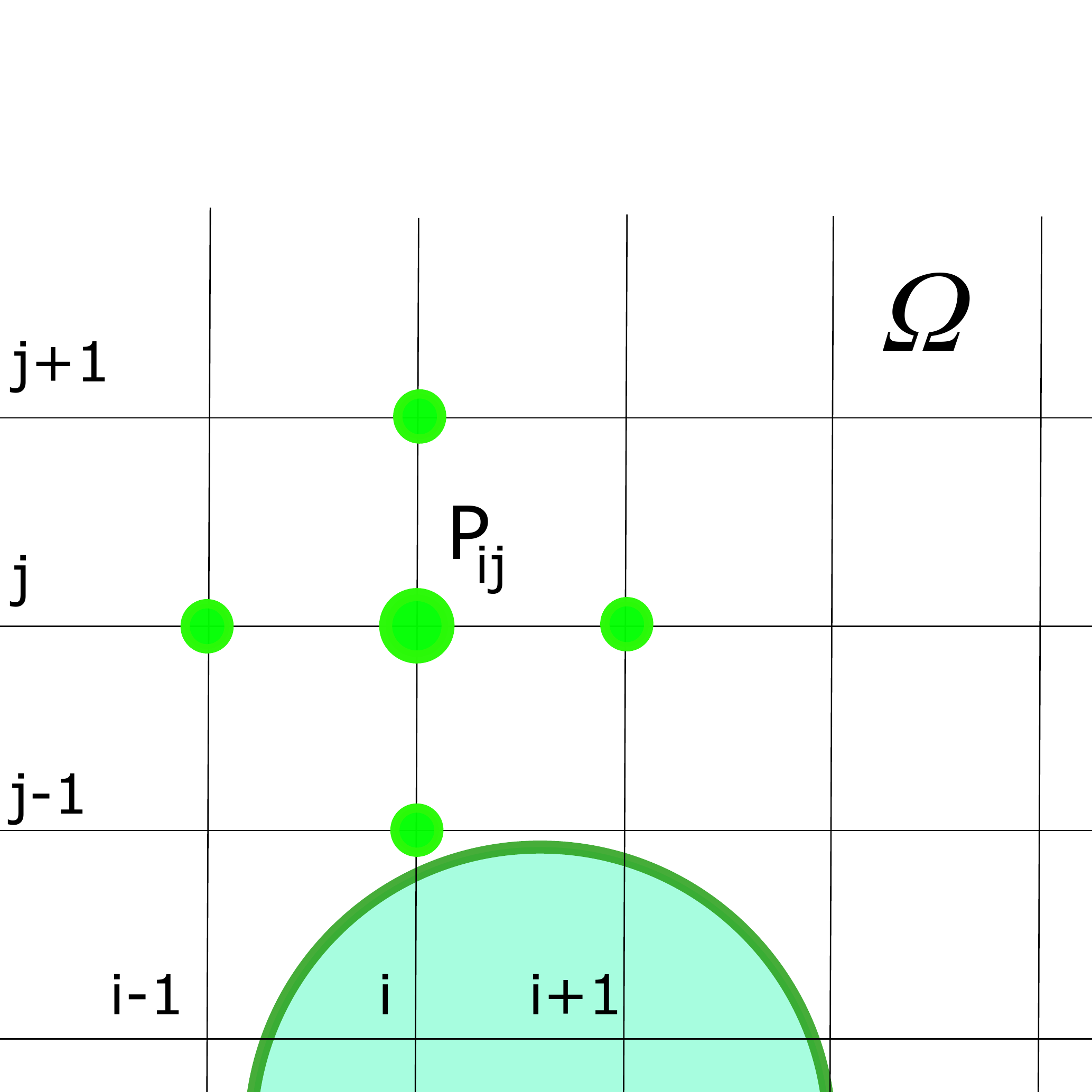}
	\end{minipage}\hfill
	\begin{minipage}
		{.45\textwidth}
		\centering
		\includegraphics[width=0.85\textwidth]{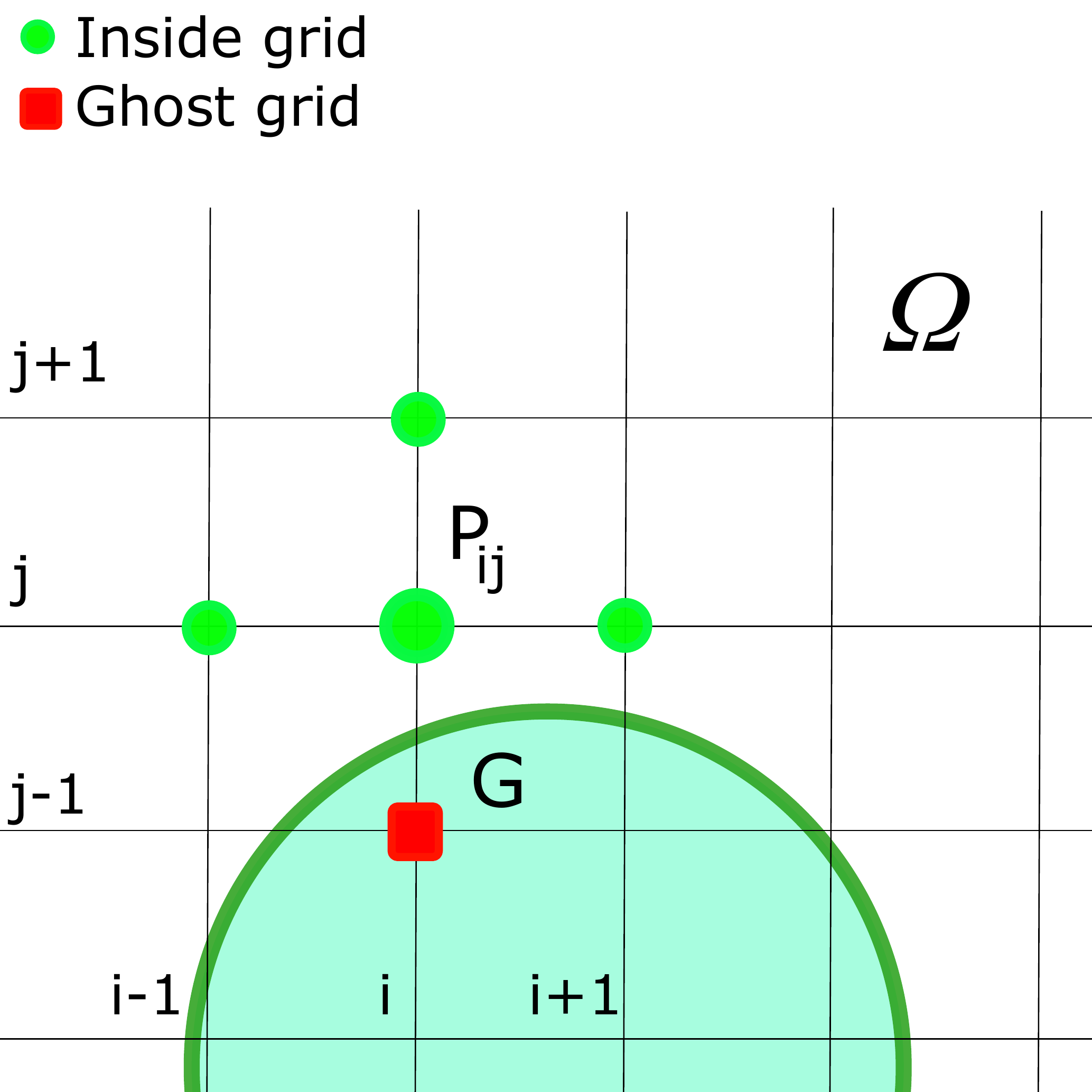}
	\end{minipage}
	\hspace*{\fill}
	\caption{\textit{Representation of the five-point stencil for the discretization of internal points $P_{ij} = (x_i,y_j)$ (left panel) and the reduced stencil when $P_{ij}$ is close to the wall $\Gamma_{\mathcal{S}}$ (right panel). In the latter case the stencil is composed by four points.}}
	\label{stencil0}
\end{figure}

If $P_{ij} = (x_i,y_j) \in \Omega^0_h$ is an internal grid point (as in Fig.~\ref{stencil0}, left panel), 
then the equation of the linear system is obtained from the discretization of the internal Eq.~\eqref{pde2d} 
and the standard central difference on a five-point stencil is used to discretize the Laplace operator on $(x_i,y_j)$. 
If $G=(x_i,y_j) \in \mathcal{G}_h$ is a ghost point, then we discretize the boundary condition \eqref{bcbubble},
following a ghost-point approach similar to the one proposed in \cite{COCO2013464}.
\begin{figure}[H]
	\centering	
	\begin{minipage}
		{.4\textwidth}
		\centering
		\includegraphics[width=\textwidth]{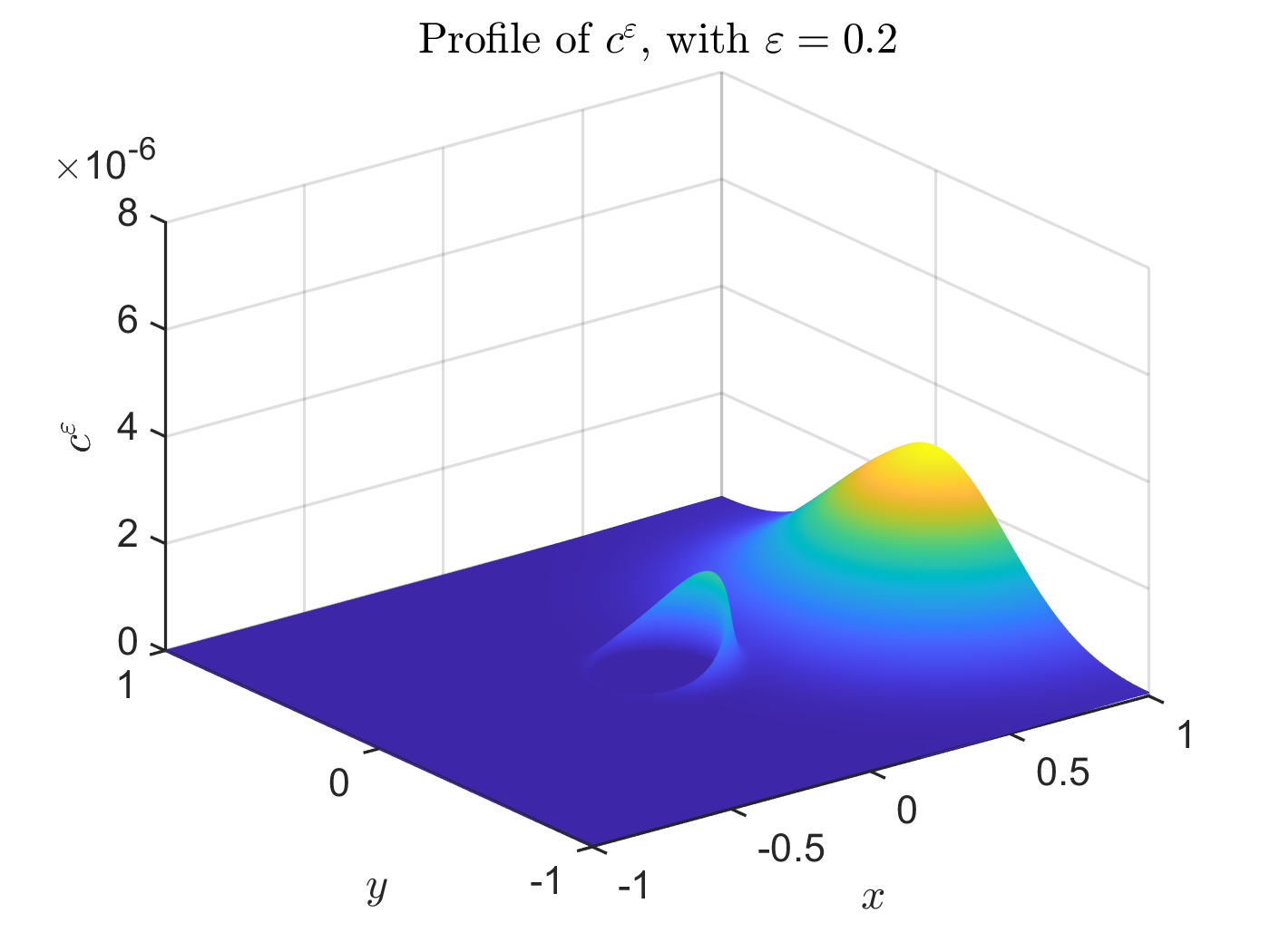}
	\end{minipage}
	\begin{minipage}
		{.4\textwidth}
		\centering
		\includegraphics[width=\textwidth]{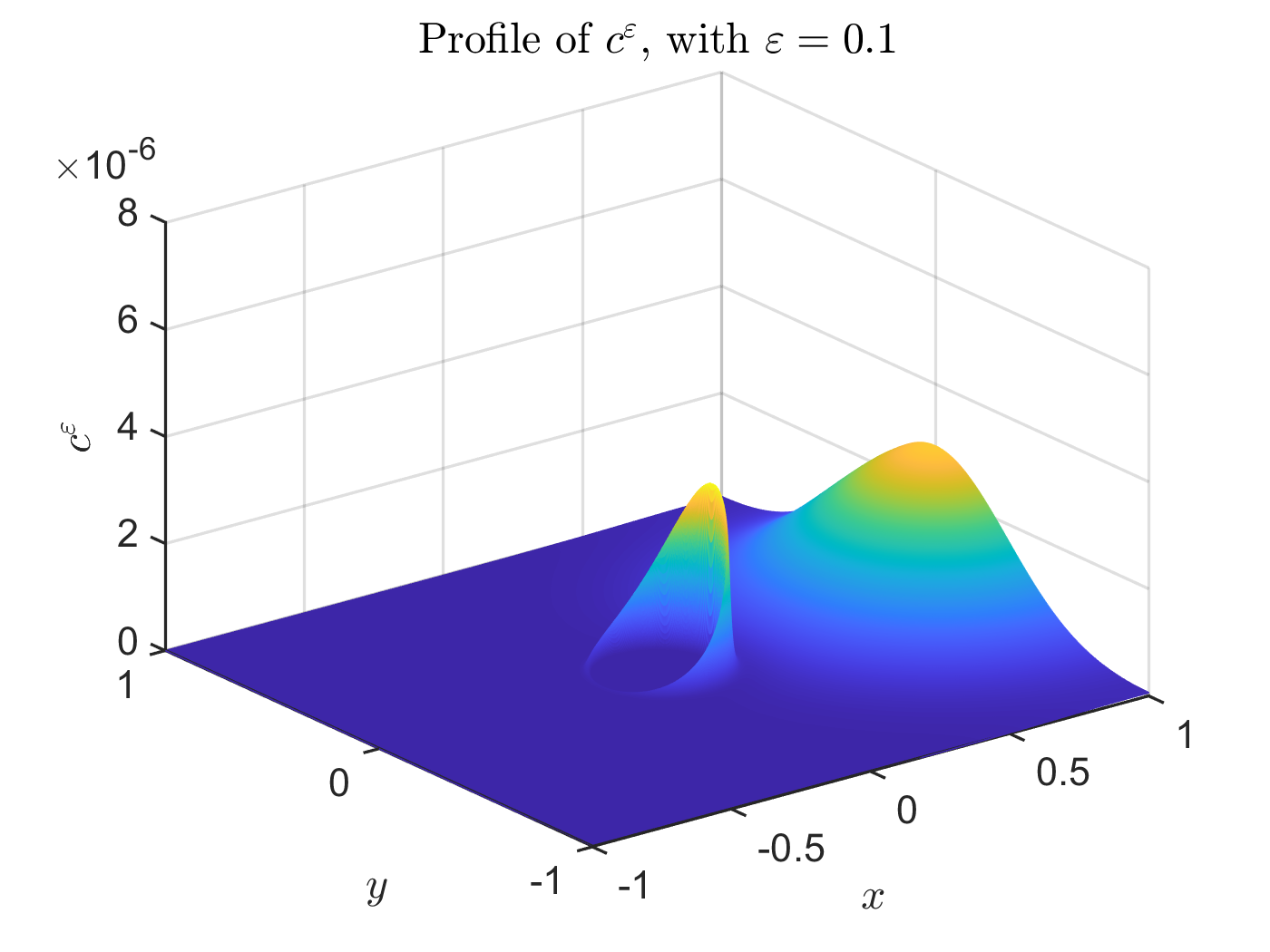}
	\end{minipage}
	\begin{minipage}
		{.4\textwidth}
		\centering
		\includegraphics[width=\textwidth]{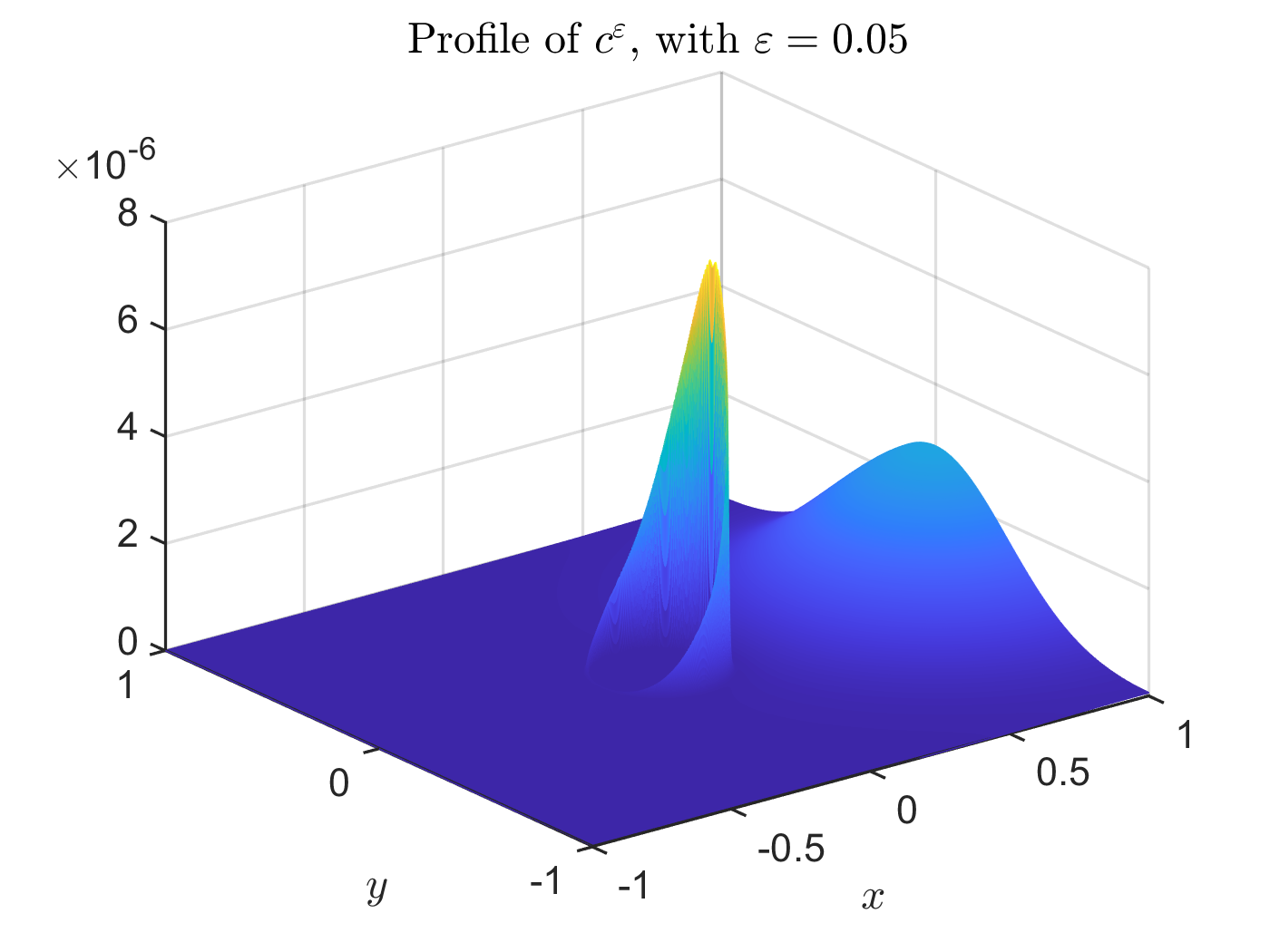}
	\end{minipage}
	\begin{minipage}
		{.4\textwidth}
		\centering
		\includegraphics[width=\textwidth]{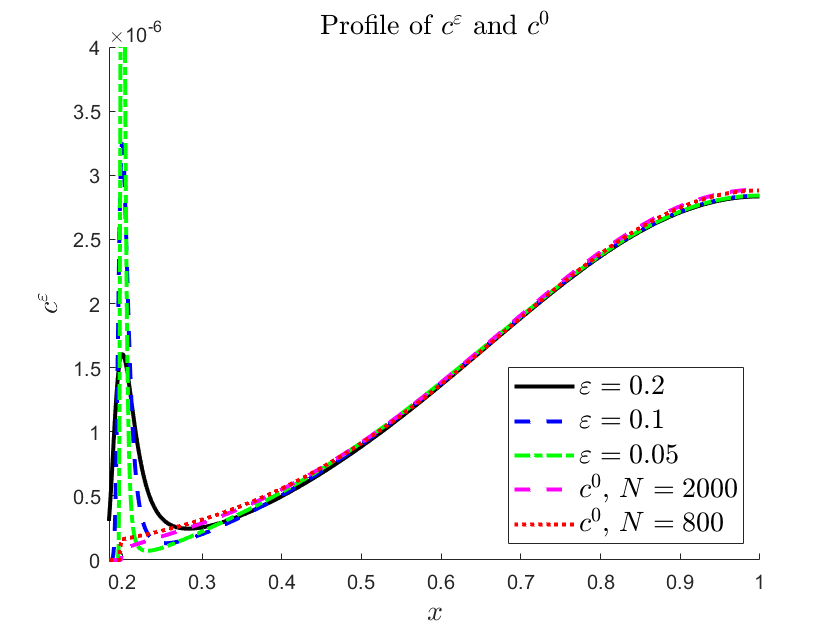}
	\end{minipage}
	\caption{\textit{In these plots we see the effect of the attractive potential in a more appropriate 2D case. Here the bubble has a reasonable circular shape. As before, in the first three panels (a,b,c) we show the solution $c^\varepsilon$ for three different values of $\varepsilon$. In panel (d) we have a 1D plot, that is a section of the previous ones for $y=0.0$ with two additional lines. The dashed black line represents the solution of the multiscale model, $c^0$, with the same space resolution of the other three cases. The dashed red line represents $c^0$ with a larger space step. The agreement of the two models improves with $\varepsilon \to 0$. Grid resolution:
			$\Delta x = 0.001$ and $\Delta x = 0.0025$ for the dashed red line, $t_{final} = 0.05$ and $M=3$. The initial condition is defined in Eq.~\eqref{ic_2d} with $x_m = 1$, $y_m = 0, \sigma = 0.1$ and $v_0 = 10^{-6}$.}}
\end{figure}

\section{Modeling the saturation effect}
\label{section_saturation}
In this section we describe how to extend the model to include nonlinear effects due to the saturation when the concentration is not negligibly small. Such effect may be relevant near the trap, because that is where the surfactant concentrates.

\subsection{Saturation: 1D model}
We derive a boundary-problem from the drift diffusion non linear equation using the same asymptotic analysis as in the case of a single species multiscale model. Starting with a one-dimensional model we impose that the coefficient of $O(\varepsilon^{-1})$ in the expansion of the 
flux is zero:
\begin{equation}
	\frac{\partial  c (x)}{\partial x}+\frac{1}{k_B T} V'(x) c (x)(1- c (x)) = 0
	\label{non_linear}
\end{equation}
This equation can be integrated giving 
\begin{equation}
	c (x) = \frac{ c_B \tilde{f}(x)}{1- c_B (1-\tilde{f}(x))}
\end{equation}
where $ c_B = c(\varepsilon L)$ is the concentration at $x=\varepsilon L$, i.e.\ at the right boundary of the support of the potential and
\begin{equation}
	\tilde{f}(x) = \exp\left(\frac{1}{k_B T}(V(\varepsilon L)-V(x))\right)
\end{equation}
denotes the classical Boltzmann factor.
The expression for the entrapped ion surface concentration $C_B $ is:
\begin{equation}
	C_B ( c_B ) = \int_{-\varepsilon}^{\varepsilon L} c (x)\,dx  = M( c_B )  c_B 
\end{equation}
where  now the coefficient $M$ in front of $ c_B $ depends on the concentration $ c_B $ just out of the potential well:
\begin{align}
	\label{eq_bc_saturation}
	M( c_B ) & = \varepsilon \int_{0}^{L+1} \frac{f(\xi,\phi)}{1- c_B (1-f(\xi))}\, d\xi, \quad f(\xi)\equiv \tilde{f}(x) &\\ 
	f(\xi,\phi) & = \exp\left( -U(\xi;\phi)\right) &
\end{align}
where we assume $U(L+1) = 0$, and
\begin{equation}
	\mathcal{I}_L(\phi,c_B)\equiv \frac{M( c_B )}{\varepsilon} = \int_{0}^{L+1}\frac{f(\xi,\phi)}{1- c_B + c_B f(\xi,\phi)}\, d\xi
	\label{I_new}
\end{equation}
{Notice that now we cannot choose $M$ arbitrarily, since 
	the integrand in \eqref{I_new} is now bounded by $1/c_B$. 
	In particular}
\[
\mathcal{I}_L(\phi,c_B) < \frac{L+1}{c_B}.
\]
This means that when taking into account the effect of saturation, two potentials with the same value of $M(0)$ will have different behavior, since a potential with a wider and less deep well can accommodate more particles, and therefore the effect of saturation will be less evident than in the case of a potential with the deeper and narrower well. 

\begin{figure}[htp]
	\centering
	\hfill
	\begin{minipage}
		{.45\textwidth}
		\centering
		\includegraphics[width=\textwidth]{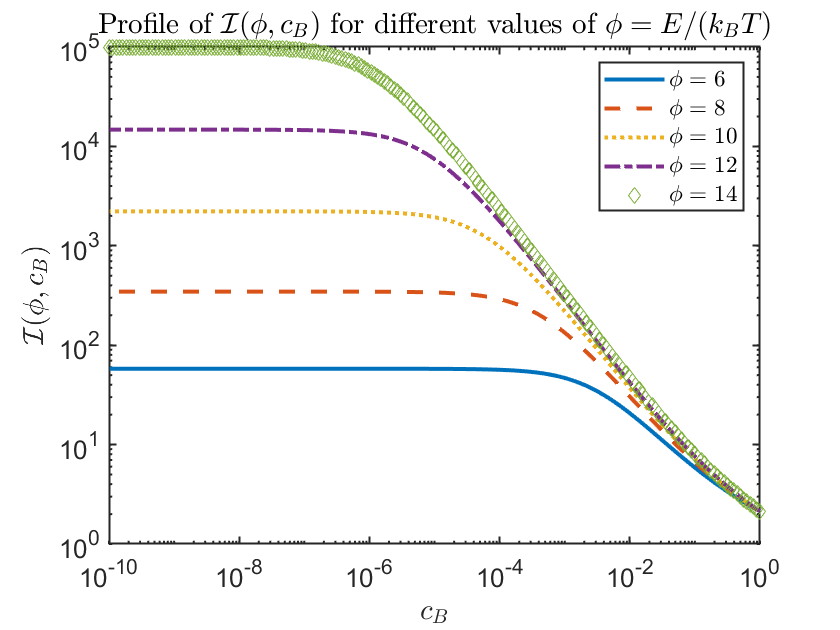}
	\end{minipage}\hfill
	\begin{minipage}
		{.45\textwidth}
		\centering
		\includegraphics[width=\textwidth]{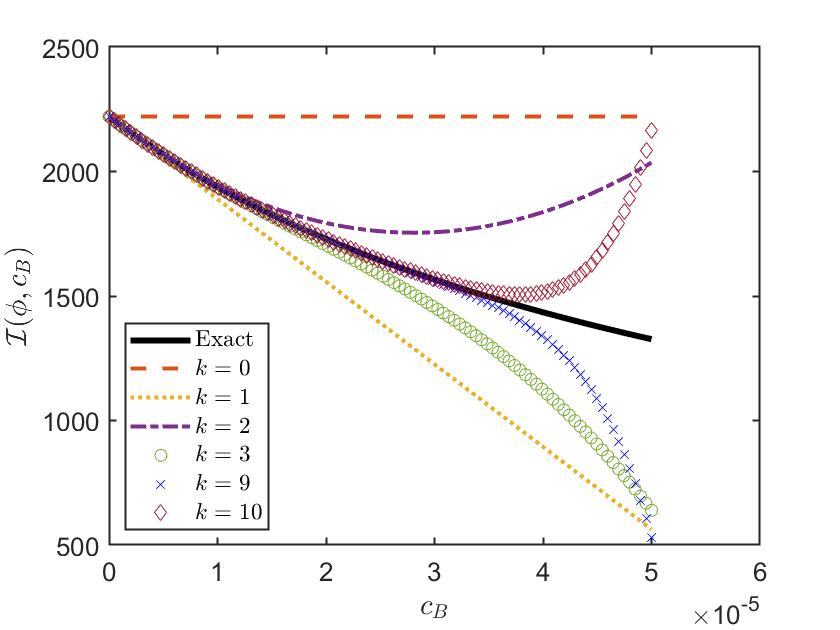}
	\end{minipage}
	\hspace*{\fill}
	\caption{\textit{On the left we plot  $\mathcal{I}(\phi,c_B)$ for different values of $\phi$. If $M$ is fixed, when $\phi$ increases, $\varepsilon$ decreases. The result is that $\mathcal{I}(\phi,c_B)$ is constant when $c_B\to 0$ and it increases {  as $\varepsilon$ decreases}. On the right we plot $M(c_B)/\varepsilon$ for $\phi = 10$ (thick black line), together with the Taylor expansion \eqref{eq:Taylor}, truncated at various order $k$. Also in Table~\ref{table_Mk} we show the oscillating behaviour of these quantities.  {The value $L=4$ has been used here.}}}
	\label{fig:Iphi}
\end{figure}

The quantity $\mathcal{I}(\phi,c_B)$ approximately represents the average amplification of the concentration in the potential well. 
In Fig.~\ref{fig:Iphi} (left) we plot the amplification $\mathcal{I}(\phi,c_B)$ for a wide range of concentrations. 

For low densities, this factor increases a lot with the parameter $\phi = E/(k_B T)$, while when the density increases, the amplification is reduced, and becomes the width of the well as $c_B$ approaches 1.

For moderate but non negligible values of $ c_B $ a truncated power series expansion of $M( c_B)$ could be used:
\begin{equation}
	M( c_B ) = \varepsilon\sum_{k=0}^\infty( c_B )^k M_k,
	\label{eq:Taylor}
\end{equation}
with
\begin{equation}
	M_k = \int_{0}^{L+1}\, f(\xi) (1-f(\xi))^k\, d\xi.
	\label{discrete_Mk}
\end{equation} 
However such an expression will not be particularly useful because the coefficients $M_k$ increase very rapidly with $k$. In Table~\ref{table_Mk} we see a few values of $M_k$ for $k = 0,\ldots,4$.
In Fig.~\eqref{fig:Iphi} (right panel) we show the expression of $M(c_B)/\varepsilon$ for $\phi = 10$, together with its Taylor expansion for several values of the order $k$ of the Taylor expansion. 
It is evident that such an expansion is not useful for concentrations bigger that about $3\times 10^{-5}$. 
\begin{table}[hbt] 
	\begin{center}
		\begin{tabular}{|c|c|c|c|c|c|} 
			\hline
			$k$ & 0 & 1 & 2 & 3 & 4 \\ \hline
			$M_k$ & 2.9065e+04  & -7.0352e+09  &  1.9534e+15 &  -5.7704e+20  &  1.7634e+26  \\  \hline 
		\end{tabular}
	\end{center}
	\caption{\textit{The first few coefficients $M_k$ (Eq.~\eqref{discrete_Mk}) of the power series  \eqref{eq:Taylor}.  $\varepsilon = 10^{-4}$.}}
	\label{table_Mk}
\end{table}

\begin{figure}[htp]
	\centering
	\hfill
	\begin{minipage}
		{.45\textwidth}
		\centering
		\includegraphics[width=\textwidth]{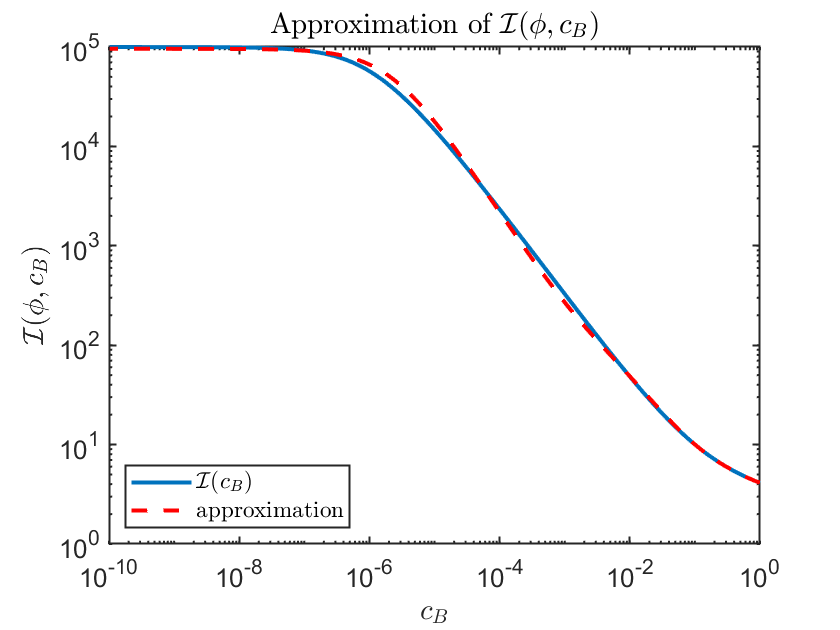}
	\end{minipage}\hfill
	\begin{minipage}
		{.45\textwidth}
		\centering
		\includegraphics[width=\textwidth]{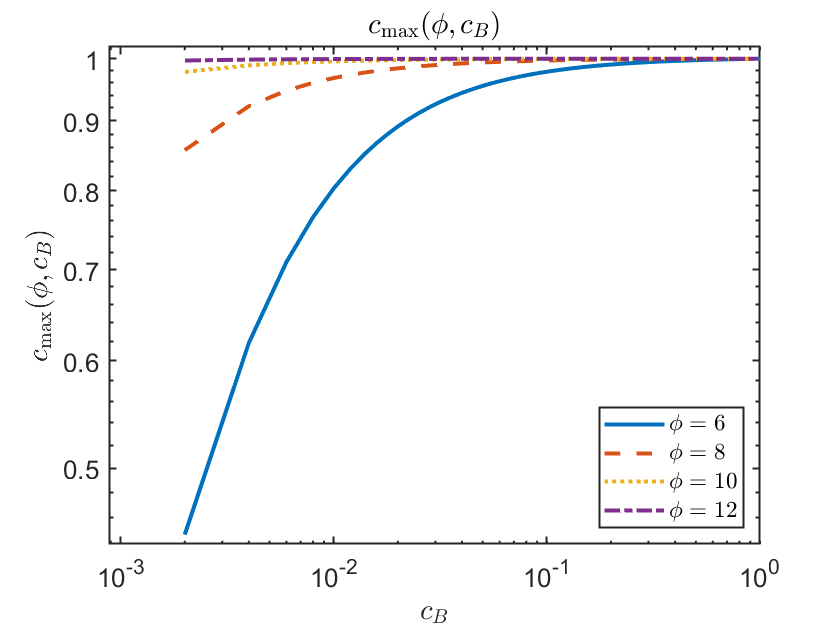}
	\end{minipage}
	\hspace*{\fill}
	\caption{\textit{On the left we plot $M(c_B)/\varepsilon$ for $\phi = 10$ (blue line), together with the rational approximation $R_{3,4}(c) = P_3(c)/Q_4(c)$ (dashed red line). In order to have a better approximation we can use rational fractions with higher degree. On the right we show the values of $c_{\rm max}(c_B)$ for $\phi \in \{6,8,10,12,14\}$. If $M$ is fixed and we increase $\phi$, the value of $\varepsilon$ goes to 0 and the maximum value of the concentration, $c_{\rm max}(\phi,c_B)$, goes very close to 1.}}
	\label{fig:Poly_approx}
\end{figure}

A more useful approximation is given by a rational function where numerator and denominator are both polynomials. In Fig.~\ref{fig:Poly_approx} (left) we see the comparison between the exact value of $\mathcal{I}(\phi,c_B)$ and a rational approximation. Of course better approximation is obtained using rational  {functions of} higher degree.

Once the trapped charge $C_B $ is expressed in terms of $c_B $ near the boundary, the boundary condition for the multiscale model with non negligible saturation effect  becomes:
\begin{equation}	
	\displaystyle \frac{d C_B }{dt} + \left. D \frac{\partial c }{\partial n}\right|_{x=0} = 0
	\label{cB_equation}
\end{equation}
where $C_B = M(c_B) c_B$.
\subsection{Results}
In this section we show various tests to validate the saturation effect described in Eq.~\eqref{eq_bc_saturation}. 

In Fig.~\ref{fig:Iphi} (left panel) we show how the quantity $\mathcal{I}_L(\phi,c_B)$ depends on $c_B$ for different values of $\phi$. We can easily see that it does not depend on $c_B$  when $c_B$ is sufficiently low and also it becomes stronger and stronger when we increase the value of $\phi$. 

\begin{figure}[tb]
		\centering
		\includegraphics[width=0.5\textwidth]{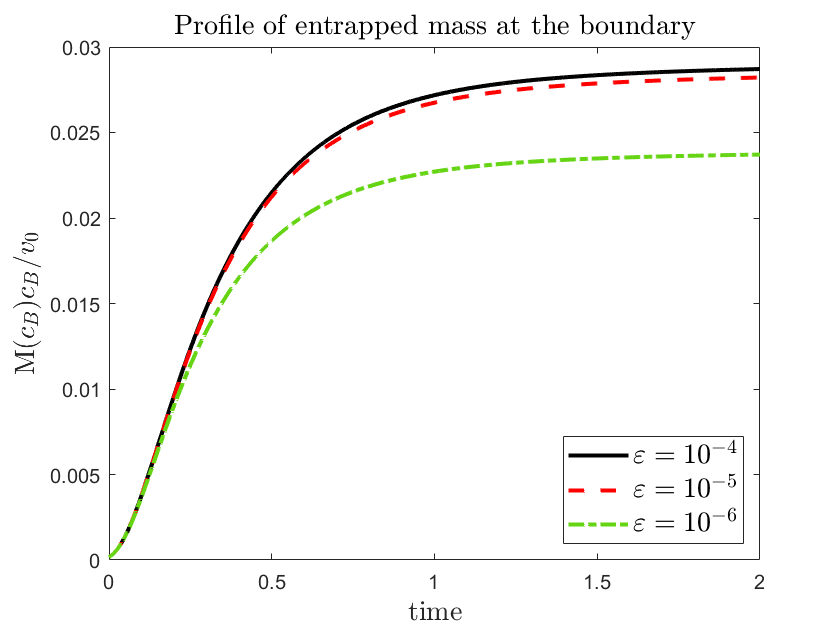}
	\caption{\textit{Left: Here we show the fraction of entrapped mass versus time for different values of $\varepsilon$. We see the solution saturates at earlier times when $\varepsilon$ decreases. Grid resolution chosen: $\Delta x = 10^{-4}$ and for the initial condition defined in Eq.~\eqref{ic_1d} we choose $v_0 = 10^{-6}, \sigma = 0.2$ and $x_m = 0.5$.}}
	\label{fig_fraction_mass}
\end{figure}

{The product $c_{\max}(\phi,c_B)\equiv c_B \mathcal{I}_L(\phi,c_B)$ provides the concentration near the bubble.}
In Fig.~\ref{fig:Poly_approx} (right panel) we show how such concentration  reaches 1 very easily when $\phi$ increases. We plot the quantity $c_B\mathcal{I}_L(\phi,c_B)$, and, as expected, we see the saturation effect becomes more effective when $\phi$ increases.
In Fig.~\ref{fig_fraction_mass} we see the fraction of entrapped mass for the Eqs.\eqref{cB_equation} versus time, for different values of $\varepsilon$. In this figure we see that the concentration saturates earlier when $\varepsilon$ decreases. 

\begin{figure}[H]
	\centering
	\hfill
	\begin{minipage}[b]
		{.33\textwidth}
		\centering
		\includegraphics[width=\textwidth]{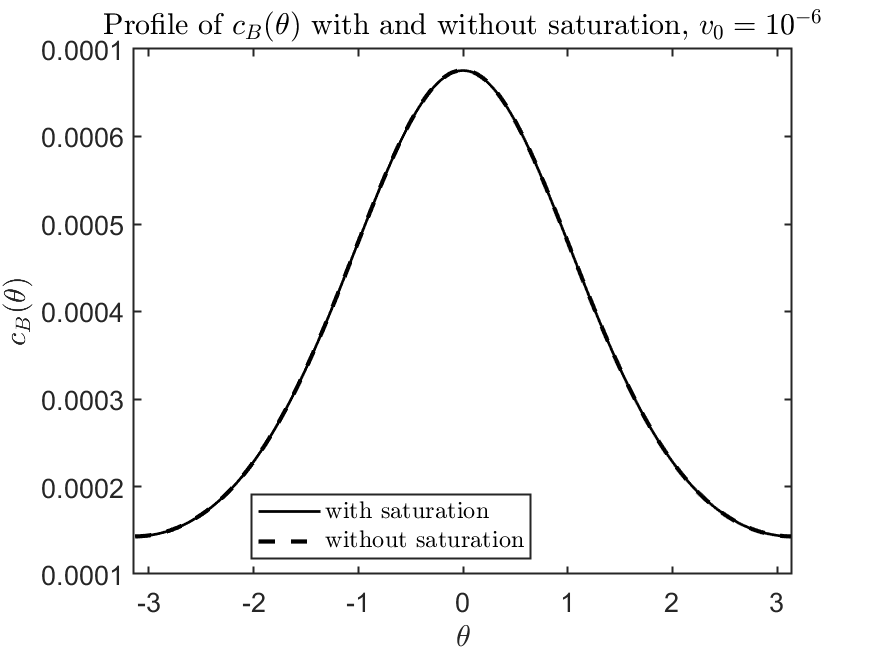}
	\end{minipage}\hfill
	\begin{minipage}[b]
		{.33\textwidth}
		\centering
		\includegraphics[width=\textwidth]{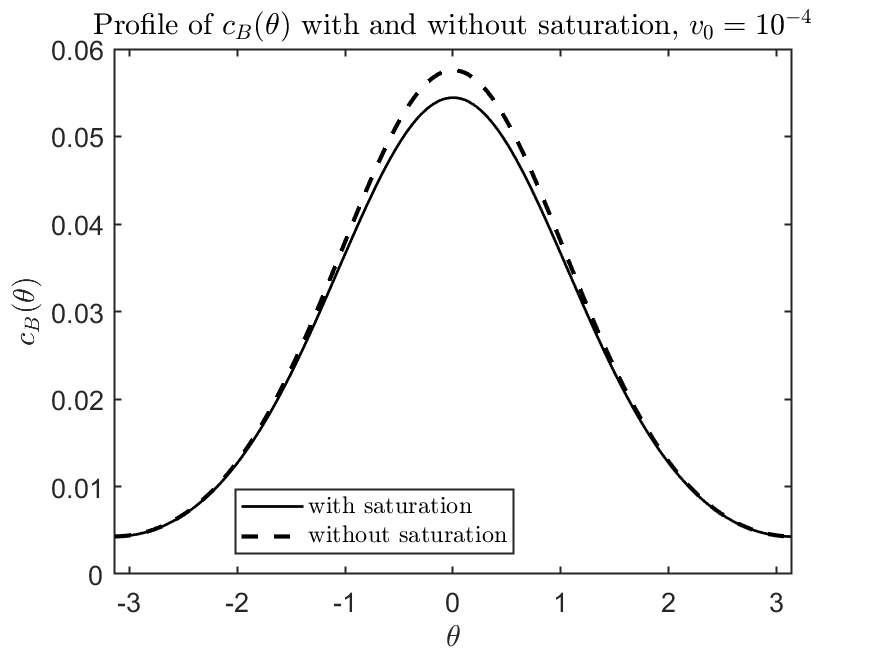}
	\end{minipage}\hfill
	\begin{minipage}[b]
		{.33\textwidth}
		\centering
		\includegraphics[width=\textwidth]{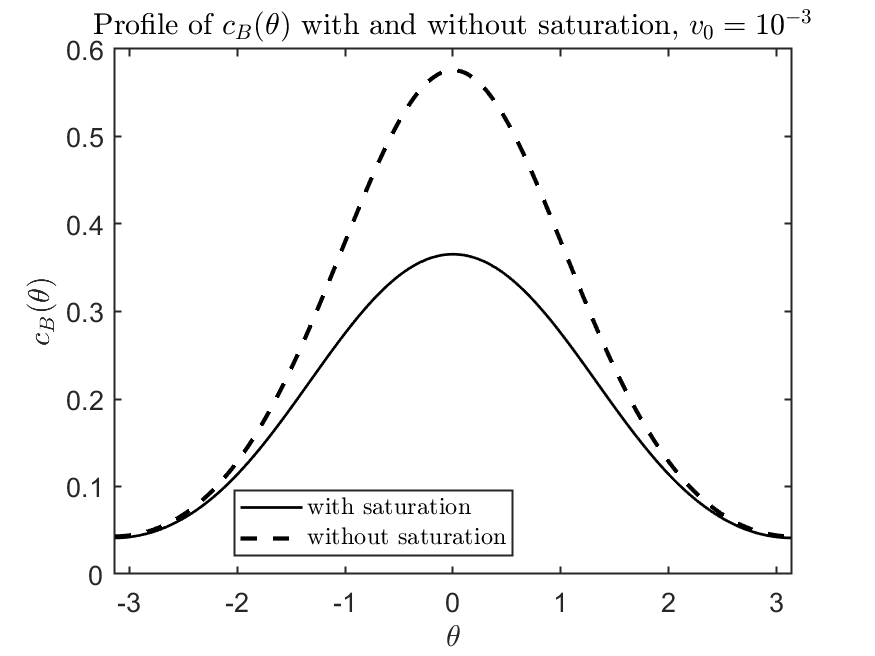}
	\end{minipage}
	\hspace*{\fill}   
	\caption{\textit{Profiles of the solution at the circular boundary, as function of the angle $\theta$. We see the saturation effect is almost negligible for very low concentrations, but it becomes more evident increasing the initial volume, $v_0$. In these plots we show the results of the solution with (black line) and without (dashed line) saturation. Grid resolution: $\Delta x = 0.002$, $t_{final} = 0.15$ and $M(c_B=0) = 3$. The values of the initial volume are described in each plot.}}
	\label{figures_cB_theta}
\end{figure}

\section{Conclusions}
\label{sec:conclusions}
In this paper we present a multiscale model for the description of the interaction between a reversible trap and a surfactant. 
The starting point is a drift-diffusion model for a single carrier, in which particles are driven by the gradient of the potential
describing the effect of the bubble.
The multiscale model is based on the assumption that the range of the potential is much smaller than typical macroscopic length scales, 
such as, for example, the radius of a spherical trap. It is shown that near the trap the particles are in local equilibrium, and therefore
the concentration satisfies a Boltzmann-type distribution. This allows a drastic simplification of  the system, which is then described 
by just a diffusion equation, with a peculiar evolutionary time dependent solution derived by the continuity equation of the trapped 
charge and the relation of local equilibrium,  {and which does not require the resolution of the small scales of the potential.}

The multiscale model solved on a uniform grid provides a tremendous reduction in the computational complexity compared to the full model, even if the latter is solved on using Adaptive Mesh Refinement, as is clearly illustrated in Appendix~\ref{A:AMR}.

Saturation effect is also taken into account.
The model is carefully validated numerically against a detailed numerical solution of a fully resolved model with a potential of width $\varepsilon$.
It is shown that the agreement between the multiscale model and the original drift-diffusion model becomes better and better as $\varepsilon$ is decreased. 
As basic potential we adopt Lennard-Jones type potential, which is uniquely determined once one specifies the depth of the potential well, that we call $E$, and its width $\varepsilon$. We observe that for very dilute solutions the effect of the potential is summarized by a single quantity that is a function of $E/k_B T$ and $\varepsilon$, because the saturation effect is negligible (for a fixed $\varepsilon$ and suitably small concentration). However we show that the effect of saturation may become relevant even for moderate concentration, since near the trap the concentration increases by orders of magnitude. This means that one has to provide realistic values of $\varepsilon$ (that we call $\vepst$ to distinguish it from the mathematical parameter used in the asymptotic expansion of the solution) and $E$, and the model is able to capture the behavior of the diffusant {\em without resolving the small scale $\vepst$}. 
This multiscale model is applied to study the evolution of the surfactant in presence of a time dependent bubble \cite{astuto2023finite}.
A physically realistic model should be based on the transport equation for both cations and ions, coupled by a self-consistent Poisson equation for the Coulomb potential \cite{CiCP-31-707}.
The current multiscale model could be used to describe carrier evolution under the condition of the so-called {\em quasi-neutral limit} \cite{CiCP-31-707,jungel}, in which the Coulomb interaction is so strong that the two carrier charge densities are practically overlapped, and move as a single carrier with effective mobility and diffusion coefficient. 
A multiscale two-carrier model which takes into account Coulomb interaction far from the quasi-neutral limit is current under investigation.

\appendix

\section{Comparison with Adaptive Mesh Refinement} 
\label{A:AMR}
{ 
	In the derivation of the multiscale model we perform an expansion in the parameter $\varepsilon$ which in our paper represents the non-dimensional thickness of the potential well near the bubble, and is given by the ratio of the actual thickness of the potential well and a typical size of the device. We can imagine that the device is of the order of a few millimeters, while a reasonable thickness of the potential well is of the order of nanometers, this giving a non-dimensional thickness of about $\vepst\approx10^{-6}$. 
	Following the argument in \cite{semplice2016adaptive}, in this appendix we give a rough estimate of the number of degrees of freedom which is necessary to fully resolve all the scales of the problem, comparing two different approaches: one based on the original drift-diffusion model in which a small grid is used only near the bubble to resolve the small scale of the thin potential of the bubble, while a coarse grid is adopted  in the bulk, where the evolution is governed by a diffusion equation, and the other one which is based on the multiscale model, in which the effect of the potential is modeled by the evolutionary boundary condition for the trapped charge,  {and a uniform coarse grid is adopted in the domain}. 
	We consider the problem in 1D, 2D and 3D. 
	Let us assume that the bulk is a domain discretized by a grid of (non-dimensional) mesh size $H$, while the region influenced by the potential has to be discretized by a grid of meshsize $h$. Let us denote by $N=1/H$ the number of grid points per direction needed to fully resolve the bulk, and $n=\vepst/h$ the number of grid points required to resolve the thickness of the potential.
	For simplicity we assume the domain is a unit cube in dimension $d$, i.e.~$\Omega = [-0.5,0.5]^d$, and the bubble is a $d$-sphere of radius $R<0.5$. 
	\begin{enumerate}
		\item In one space dimension one needs $n$ grid points to discretize the region influenced by the potential, and $N$ grid points to discretize the bulk.
		\item In two space dimensions we assume we have small regions of size $h^2$ inside the potential and large regions of size $H^2$ in the bulk. The number of regions inside the potential is then given by the ratio between the area influenced by  the potential and the area of the small cell: $2\pi R\vepst/h^2$. On the other hand, the number of elements needed to model the boundary conditions in the multiscale model is given by $2\pi R/H = 2\pi R N$.
		The number of elements in the bulk is approximately $N^2$ in both cases.
		\item In three space dimensions the number of elements in the region influenced by the potential is given by the volume of such region, $4\pi R^2 \vepst$, divided by the volume $h^3$ of each small element. 
	\end{enumerate}
	We summarize in Table~\ref{tab:compare} the number of degrees of freedom which are needed to fully resolve the problem in the two approaches (AMR and Multiscale). Here we neglected the effect of factors $2\pi R$ and $4\pi R^2$ which are $O(1)$.}

\begin{table}[]
	\centering
	\begin{tabular}{|c|c|c|c|c|} \hline
		d  & \multicolumn{2}{|c|}{AMR \#dof}  & \multicolumn{2}{|c|}{Multiscale \#dof}         \\ \hline
		1  & $n+N$               & $110$     & $1 + N$      & $101$      \\ \hline
		2  & $n^2/\vepst + N^2$  & $10^8$    & $N + N^2$    & $10^4$     \\ \hline
		3  & $n^3\vepst^2 + N^2$ & $10^{15}$ & $N^2 + N^3$  & $10^6$     \\ \hline
	\end{tabular}
	{ \caption{Order of magnitude of the number of degrees of freedom needed to resolve the problem in dimension $d$, using either full model discretized by adaptive mesh refinement (AMR) or discretized by the new multiscale model. The number in columns 3 and 5 refer to the specific assumption $n=10$, $N=100$ and $\vepst = 10^{-6}$.}
		\label{tab:compare}
	}
\end{table}
{	From these considerations it follows that the multiscale model allows an improvement of efficiency by orders of magnitude over a detailed direct computation performed by resolving all the small scales of the problem. }

\section{Laplace-Beltrami operator}
\label{A:LB}
{ 
	
	Here we prove that 
	\begin{equation}
		\label{A:LapBel}
		\partial_i\tpar_i c = \tpar_i\tpar_i c =: \Delta_\perp c
	\end{equation}
	where $\Delta_\perp$ is the Laplace-Beltrami operator. We also provide a representation of $\Delta_\perp$ in Cartesian coordinates.
	
	First, observe that
	\[
	\partial_i = \tpar_i  + n_in_j\partial_j,
	\]
	therefore 
	\[
	\partial_i\tpar_i c = \tpar_i\tpar_i c + n_in_j\partial_j\tpar_i c.
	\]
	Now consider the term
	\[
	n_in_j\partial_j\tpar_i c = n_i\pad{}{n}\tpar_i c.
	\]
	Since $n_i \tpar_i c \equiv 0$, it is
	\[
	0 = \pad{}{n}(n_i\tpar_i c) = \pad{n_i}{n}\tpar_i c + n_i \pad{}{n}\tpar_i c.
	\]
	Now $\partial n_i/\partial n = n_j\partial_j n_i = 0$, therefore 
	\[
	n_i\pad{}{n}\tpar_i c = 0
	\]
	from which \eqref{A:LapBel} follows.
	
	Now we give a representation for $\Delta_\perp$.    
	We start observing that $\tpar_i c = h_{ij}\partial_j c$ is orthogonal to $\Sigma_B$, since
	$n_ih_{ij} = 0$. 
	It is
	\begin{align}
		\Delta_\perp c = \tpar_i\tpar_i c &  = \tpar_i(h_{ij}\partial_j c) = \tpar_i(\partial_i c - n_i n_j \partial_j c) \\
		& = h_{ij}\partial^2_{ij} c - \tpar_in_i\pad{c}{n} - n_i (\tpar_i n_j) \partial_j c - n_i n_j \tpar_i\partial_j c
	\end{align}
	
	The last two terms are zero, because 
	$n_i\chi_{ij} = 0$ and $n_i\tpar_i\partial_j c = 0$.
	
	Therefore it follows:
	\[
	\Delta_\perp c = \tpar_i\tpar_i c  = 
	h_{ij}\partial^2_{ij}c - \chi_{ii} n_j\partial_j c
	\]
	We recall that the trace of the second fundamental form of a 2D surface in 3D is the sum of the two principal curvatures, i.e. twice the mean curvature of the surface. 
	
	\subsection*{Laplace-Beltrami operator in two dimensions}
	In two space dimensions (see Fig.\ref{fig:bubble}) expressions considerably simplify.
	Here we prove Eq.~\eqref{eq:BL2D}, i.e. 
	\begin{equation}
		\frac{\partial^2 c}{\partial \tau^2} \equiv
		\tau_j\partial_j(\tau_i\partial_i c) = 
		h_{ij}\partial^2_{ij}c - \chi_{ii}\pad{c}{n} \equiv \Delta_\perp c.
		\label{A2:LB}
	\end{equation}
	Indeed it is 
	\begin{equation}
		\label{A:tau2}
		\frac{\partial^2 c}{\partial \tau^2} = 
		\tau_j\partial_j(\tau_i\partial_i c) = 
		\tau_i\tau_j\partial^2_{ij} c + \tau_j(\partial_j\tau_i)\partial_i c = 
		\tau_i\tau_j\partial^2_{ij} c - \frac{1}{R} n_i\partial_i c
	\end{equation}
	the latter equality is a consequence of the relation 
	\[
	\pad{\hat{\tau}}{\tau} = -\frac{1}{R} \hat{n}
	\]
	where $R$ denotes the local radius of curvature.
	Finally we observe that in 2D it is 
	\begin{equation}
		\label{A:h}
		h_{ij} \equiv \delta_{ij}-n_in_j = \tau_i\tau_j
	\end{equation}
	as it can be easily shown by inspection checking that 
	\[
	n_i n_j + \tau_i \tau_j = \delta_{ij}
	\]
	Using Eqs.~\eqref{A:tau2} and \eqref{A:h} in \eqref{A2:LB} completes the proof of relation \eqref{eq:BL2D}.
}

\bibliographystyle{unsrt}
\bibliography{main}

\end{document}